\documentclass{article}
\usepackage{graphicx} 
\usepackage{enumitem}
\usepackage{comment}

\usepackage[a4paper,
            bindingoffset=0.2in,
            left=1in,
            right=1in,
            top=1in,
            bottom=1in,
            footskip=.25in]{geometry}
            
\usepackage{amsmath, amssymb}
\usepackage{tikz}
\tikzset{>=latex} 
\usepackage{pgfplots} 
\pgfplotsset{compat=1.18} 
\usetikzlibrary{patterns, intersections,backgrounds,positioning,calc}
\usepackage{xcolor}
\definecolor{color1}{HTML}{66c2a5}
\definecolor{color2}{HTML}{fc8d62}
\definecolor{color3}{HTML}{8da0cb}
\definecolor{color4}{HTML}{e78ac3}
\definecolor{color5}{HTML}{a6d854}
\definecolor{color6}{HTML}{ffd92f}
\definecolor{color7}{HTML}{e5c494}
\definecolor{color8}{HTML}{b3b3b3}
\usepackage[outline]{contour} 
\contourlength{1.2pt}

\usepackage{pgfplots}
\usepackage{pgfplotstable}
\usepackage{amsmath}
\usepackage{caption}
\usepackage{subcaption}
\pgfplotsset{compat=1.18}

\pgfmathdeclarefunction{gauss}{2}{%
  \pgfmathparse{1/(#2*sqrt(2*pi))*exp(-((x-#1)^2)/(2*#2^2))}%
}

\pgfmathdeclarefunction{gausscdf}{2}{%
  \pgfmathparse{1/(1 + exp(-0.07056*((x - #1)/#2)^3 - 1.5976*(x - #1)/#2))}%
}
\pgfplotstableread{
x y
0.05 0.1
0.20 0.2
0.35 0.3
0.50 0.4
0.60 0.5
0.70 0.6
0.80 0.7
0.90 0.8
0.95 0.9
1.00 1.0
}\invcdfdata

\usepackage{natbib}
\usepackage{eurosym}
\usepackage[font=small]{caption}
\usepackage{subcaption}

\usepackage{gensymb}

\usepackage{hyperref}

\usepackage{multicol}
\usepackage{wrapfig}

\usepackage{mathtools} 

\title{Electrolyzers Bidding in Electricity Markets under Green Hydrogen Regulations and Uncertainty}
\author{Andrea Gloppen Johnsen, Lesia Mitridati, Jalal Kazempour, Line Roald}
\date{July 2025}

\begin{document}
\let\texteuro\euro

\maketitle

\begin{abstract}
    Hydrogen produced through electrolysis offers a pathway to decarbonize hard‑to‑abate sectors by replacing \textit{gray} hydrogen derived from natural gas reforming when produced using renewable power. However, grid‑connected electrolyzers may inadvertently increase power‑system emissions, resulting in hydrogen whose life‑cycle intensity is similar to or higher than that of gray hydrogen. To address the high cost barrier of electrolytic hydrogen, both the E.U. and U.S. have introduced subsidy schemes conditional on low associated emissions. One key requirement is temporal matching, under which a subsidy applies only to the hydrogen volume that, ex-post, can be shown to match renewable generation over each one‑hour interval. This requirement exposes the electrolyzer to uncertainty in the subsidy‑eligible volume and thus the value of the produced hydrogen.
This paper develops an uncertainty‑aware day‑ahead bid curve for a grid‑connected electrolyzer. We formulate a linear program that maximizes expected profit across scenarios of renewable production and derive the bid curve from its Karush‑Kuhn‑Tucker conditions. A case study demonstrates that incorporating renewable uncertainty into the bid curve increases electrolyzer profit by approximately 4\%, although it does not improve ex-post temporal matching. This finding highlights a potential distortion in the incentive effects of temporal‑matching regulations when uncertainty is taken into account.
\end{abstract}

\section{Introduction}
\subsection{Electrolyzers enable indirect electrification in hard-to-abate sectors}
In order to reach net-zero emission targets, it is necessary to decarbonize across the entire energy system. For many sectors, decarbonization can be achieved through parallel electrification and build-out of renewable energy sources (RES). However, direct electrification is difficult or impossible in certain sectors,  such as heavy transportation and chemical manufacturing, due to high costs or technological challenges. In some of these sectors, one pathway to decarbonization is through indirect electrification, using electrolytic hydrogen.  Electrolytic hydrogen is produced by \textit{electrolyzers} through the consumption of electricity and water. 
Currently, hydrogen is produced mainly by reforming of natural gas, and  is typically referred to as \textit{gray hydrogen} with emission intensities estimated at around 10 kg of CO$_2$ per kg H$_2$ \citep{InternationalEnergyAgency2024Global2024}. Producing the current worldwide consumption of 95 Mt of gray hydrogen using existing electrolysis technology instead would take about 5,000 TWh of electricity,\footnote{Assuming 19 kg$_{\rm H_2}$ per MWh of electricity, which corresponds to a 63\%  conversion efficiency with respect to the lower heating value of hydrogen.} which is comparable to the current annual electricity consumption of the United States \citep{U.S.EnergyInformationAdministrationStateProfiles}. 

\subsection{The challenge of increases in power system emissions}
Therefore, while electrolyzers present an opportunity to decarbonize hard-to-abate sectors, they also pose a major challenge to the decarbonization of power systems, as their deployment will significantly add to the already growing electricity demand. As power systems transition from conventional to renewable generation, added electrolyzer demand risks increasing the power system emissions if not all the additional demand can be met by RESs. Studies show that grid-connected electrolyzers can result in  emission intensities more than three times that of gray hydrogen \citep{deKleijne2024WorldwideTransport}. It is evident that the emission intensity of electrolytic hydrogen depends on the emission intensity of the electrical energy consumed. However, literature shows that exact emissions accounting of grid consumption can be challenging \citep{Langer2024DoesImpact}, and there has been vigorous discussion regarding under what conditions the hydrogen from grid-connected electrolyzers should qualify as \textit{renewable} or \textit{green} hydrogen. 

To ensure a net-reduction of emissions, both the E.U. and U.S.  have formulated regulation of electrolytic hydrogen production \citep{EUCommission2023DelegatedOrigin, InternalRevenueServiceTreasury2023SectionProperty}. Similar in both regulations is that electrolytic hydrogen can qualify for a subsidy if the electrolyzer's electricity consumption complies with requirements on \textit{additionality}, \textit{locational matching}, and \textit{temporal matching}. In general, this means that the electrolyzer must prove that its consumption \textit{matches} with the production of one or several RES, through power purchase agreements (PPAs), guarantees of origin (GOs), or renewable energy certificates (RECs). The additionality requirement refers to that the matching RES must be newly built as to avoid cannibalization of the existing RESs. The locational matching requirement refers to that the electrolyzer and matching RES have geographical correlation, e.g. that they are located in the same bidding zone. Finally, the temporal matching requirement refers to that the electrolyzer consumption must be temporally correlated to the matching RES. However, due to the uncertain and variable nature of RES, temporal matching cannot easily be enforced in real time. Instead, different time-frames have been proposed, starting from yearly or monthly matching and transitioning to hourly matching. 

\subsection{Hourly matching avoids power system emissions increases}
Hourly matching has been shown in literature to efficiently mitigate increases in power system emissions as a result of the additional load from electrolyzers, whereas monthly and annual matching have shown varying effects. \cite{Ricks2023MinimizingStates} show that annual matching consistently leads to emission intensities higher than that of gray hydrogen, while hourly matching can achieve near zero emission intensities in areas with relatively-low competition for renewable energy sources. \cite{Zeyen2024TemporalHydrogen} also show that hourly matching will avoid increases in power systems emissions due to the additional load of electrolyzer. However, they further find that depending on the generation mix in the power system and the demand-side flexibility of the electrolyzer, monthly or annual matching could also avoid large increases in emissions, keeping the emission intensity of the electrolytic hydrogen below that of gray hydrogen. \cite{Giovanniello2024TheProduction} show that how the requirement of \textit{additionality} is interpreted and modeled can have a great effect on the resulting emission intensities of electrolytic hydrogen under an annual time matching requirement, while hourly matching requirements consistently keep the emission intensity of electrolytic hydrogen below that of gray hydrogen. \cite{Langer2024DoesImpact} studies the emission effects of various REC criteria in a review article (including, i.a., the works by \cite{Ricks2023MinimizingStates}, \cite{Zeyen2024TemporalHydrogen} and \cite{Giovanniello2024TheProduction}), in a general study for any electrical load. They find that annual matching cannot guarantee a reduction in power system emissions due to the cannibalization of RES projects (i.e., that the load purchases RECs from RES that would have been built in either case). However, they too find that hourly matching significantly reduces emissions across the reviewed studies.  


\subsection{The cost of hourly matching }
While hourly matching shows promising results on avoiding increases in electrolyzer induced emissions, it has been opposed by industry due to concerns on its effect on the business case of electrolyzers \citep{FertilizersEurope2025CallProduction}. As RES production is typically variable, the electrolyzer might have to reduce its average utilization factor in order to comply with the hourly matching requirement. There is a concern that since electrolyzers are capital-intensive, reducing the overall hydrogen production output would increase the levelized cost of hydrogen (LCOH), thus reducing their cost-competitiveness.   The literature discuss to what extent time-matching requirements will increase the production costs of electrolytic hydrogen. 
\cite{Ricks2023MinimizingStates} report that hourly matching can be achieved while increasing the LCOH with only \$1 per kg compared to unconstrained electrolyzer operation. Similarly, \cite{Ruhnau2023FlexibleEmissions} finds that the LCOH would increase from EUR3.3 per kg to EUR4.5 per kg when hourly matching is enforced. \cite{Zeyen2024TemporalHydrogen} find that hourly matching would only increase costs with 7-8\% compared to monthly matching as long as cheap hydrogen storage technologies are available, i.e., under demand-side flexibility.  

\subsection{Research gaps}

While \cite{Ruhnau2023FlexibleEmissions} takes the perspective of an investor, and \cite{Ricks2023MinimizingStates} and \cite{Zeyen2024TemporalHydrogen} take the perspective of a system operator, all papers model the hourly-matching requirement as a constraint in an optimization model, enforcing that the electrolyzer never consumes more (or less) than the
Common in all reviewed literature is that they base their conclusions on models that assume perfect information of the matching RES production of every hour. However, the actual production from a RES will only be known ex-post. The electrolyzer must therefore plan its consumption schedule based on forecasted RES production, which, depending on the availability of updated information, market commitments, and flexibility in hydrogen demand, might prove difficult to differ from in real time. Therefore, although the electrolyzer \textit{intends} to match its consumption with RES production, the uncertainty of the matching RES might cause disparities between the consumed and produced energy upon realization. 


\subsection{Optimal bidding under uncertainty}
Under the proposed regulations, grid connected electrolyzers that participate in electricity markets are exposed to the uncertainty of the time-matching RES, which propagates to the \textit{value of its produced hydrogen}. The consumption that can be shown ex-post to have time-matched with RES(s) production will qualify for a subsidy, which inflates the value of the produced hydrogen. If the electrolyzer consumption exceeds the production of the time-matched RES(s), this exceeding hydrogen does not qualify for the subsidy and is only valued at its market rate or a contracted price. This forms a problem of optimal bidding under uncertainty for an electrolyzer participating in electricity markets.

Optimal procurement under uncertainty represents a classic challenge in operations research, often exemplified by the newsvendor problem, which determines the order quantity that maximizes expected profit given stochastic demand and known cost parameters~\citep{Khouja1999TheResearch,Choi2012HandbookProblems}.  In this problem, wrongly estimating the demand carries a cost, either in the form of unsold inventory if demand was overestimated (i.e. remaining newspapers)  or lost revenue if the demand was underestimated, i.e., an opportunity cost. 
Applications of the newsvendor problem in electricity markets include the work by \cite{Pinson2007TradingPower} applying the newsvendor problem to a wind farm participating in the day-ahead and balancing markets. There, the ``demand” corresponds to uncertain wind production at the time of bidding; overestimating output exposes the producer to imbalance penalties, whereas underestimating output forgoes potential revenue. This supply‑side uncertainty in renewable generation mirrors the demand-side uncertainty in the classic newsvendor formulation. 

In the case of an electrolyzer participating in the day‑ahead market under temporal‑matching regulation, uncertainty arises from the availability of subsidies. Overestimating matched renewable output may lead the electrolyzer to purchase electricity at prices that prove unprofitable when the subsidy is not available. Conversely, underestimating the matched output can cause the electrolyzer to miss opportunities to buy power at sufficiently low prices, thereby reducing potential profits. 

Many electricity markets allow participants to submit \textit{bid curves}, i.e., a curve stating for what price they are willing to buy or sell a certain quantity of energy \citep{NordPoolSingleOrder}. Typically, to maintain convexity when solving the market clearing problem, these curves can be provided as pice-wise linear step-curves, with decreasing steps for consumers and increasing steps for producers \citep{NEMOCommittee2019EUPHEMIAAlgorithm}. Therefore, the electrolyzer's bid curve to the electricity market should follow this requirement, with a step-wise increasing price-quantity curve.

\subsection{Contribution}
In this work, we seek to investigate how renewable generation uncertainty impacts the profits of an electrolyzer in the context of hourly matching requirements. To enable this analysis, we derive a bid curve  for an electrolyzer that takes the potential subsidy achieved when  time-matched as well as the uncertainty of the matching RES into account. This bid curve can be used to participate in, e.g., the day-ahead electricity market. We derive the bid curve from the Karush-Kuhn-Tucker (KKT) conditions\footnote{The KKT conditions are a set of equality and inequality conditions that must be fulfilled for a given decision to be optimal in its original problem. See \cite{Boyd2024ConvexOptimization}, chapter 5.5 \textit{Optimality Conditions}.
} of an optimization problem, which maximizes the expected profit with uncertain value of hydrogen production due to the uncertainty of the time-matched RES. We apply this bid curve to a case study of Denmark, utilizing historical market data, and compare the derived uncertainty-aware bid curve to those of a deterministic and a perfectly informed curve. Our results show that there is a 9\% loss in electrolyzer profits between the perfect information curve and the deterministic curve, while the uncertainty-aware curve reduces this loss to 6-4\% under naive information about the RES distribution. However, our results also show that bidding in an uncertainty-aware manner leads to more consumption that exceeds the realized production of the matching RES compared to the deterministic case. Due to the framework of this study and the availability of data we make no investigation into how this increase in exceeding consumption impacts system emissions, but conclude that its effect cannot be ruled out based on the current study. We recommend that to properly evaluate the consequences of hourly time-matching requirements on the cost of hydrogen and emissions of the power system, the uncertainty of RES production and ways to mitigate it should be taken into account.

\section{Results}

\subsection{Case study: West-Denmark 2024}
In order to evaluate the performance of our derived bid curves, we form a case study on historical day-ahead price and wind (forecast and realization) data from the western Danish bid zone (DK1) in 2024. For every hour of the year, we derive a bid curve based on the day-ahead wind-power forecast. This bid curve is evaluated against the day-ahead electricity price of the hour, which determines the electrolyzer's operational point. This approach assumes that the electrolyzer is a price taker, i.e., that its consumption would not have a significant effect on the final clearing price. 

We assume a constant electrolyzer efficiency $\eta$ equal to 18 kg$_{\rm H_2}$ per MWh$_{\rm el}$, a market rate hydrogen price $\pi^{\rm gray}$ of EUR 2 per kg$_{\rm H_2}$ and compliance subsidy $\pi^{\rm green}$ of EUR 4 per kg$_{\rm H_2}$. We calculate the utility, i.e., the MWh$_{\rm el}$ equivalent value of gray and green hydrogen as follows:
\begin{align}
    &\lambda^{\rm gray} = \eta \pi^{\rm gray}, \label{gray_price} \\
    & \lambda^{\rm green} = \eta (\pi^{\rm gray} + \pi^{\rm green}). \label{green_price} 
\end{align}

In the following, we distinguish between the hydrogen \textit{price}, which is given per kg$_{\rm H_2}$ and noted by $\pi$, and the hydrogen \textit{value}, given per MWh$_{\rm el}$ and noted by $\lambda$. 

\begin{figure}[t!]
    \centering

    \begin{minipage}[t]{0.48\textwidth}
        \centering
        \includegraphics[width=\linewidth]{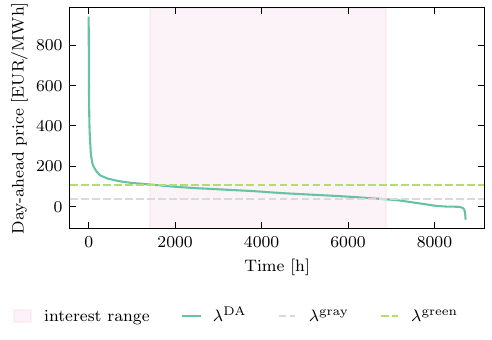}
        \caption{Price-duration curve of the considered case study of DK1 2024. $\lambda^{\rm DA}$ shows the price of every hour in the year, sorted in descending order. $\lambda^{\rm gray}$ and $\lambda^{\rm green}$ show the electrolyzer's utility when producing gray and green hydrogen respectively, while the shaded area highlights the part of the duration curve that falls within this price-range of interest.}
        \label{fig:price_duration}
    \end{minipage}
    \hfill
    \begin{minipage}[t]{0.48\textwidth}
        \centering
        \includegraphics[width=\linewidth]{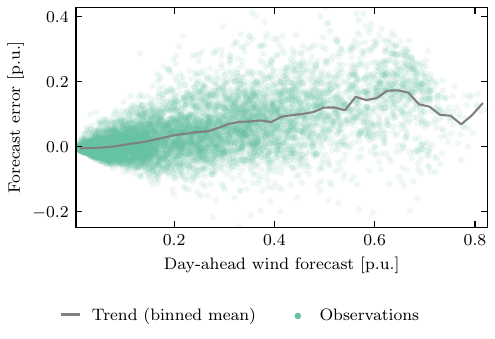} 
        \caption{Onshore wind day-ahead forecast and the associated forecast error for the considered case study of DK1 2024. There is a positive bias on the mid-to-high forecasted values, i.e., higher realization than forecast. The data is given for the aggregated onshore wind capacity, and per unit (p.u.) of the highest dataset value.}
        \label{fig:wind_distribution}
    \end{minipage}

\end{figure}

Figure \ref{fig:price_duration} shows the price duration curve of the case study. The day-ahead price in each hour of the year is plotted in descending order. We show the electrolyzer utility of producing gray and green hydrogen, marked by dashed, horizontal lines, for reference. These lie at EUR 36 and EUR 108 per MWh, respectively. As our method of deriving the electrolyzer bid curve assumes truthful bidding, i.e., bidding after the true marginal value of hydrogen production, any price that fall above the green or below the gray hydrogen value is out of interest in the study, as these represent the minimum and maximum value of hydrogen. In these cases, any derived bid curve will result in the same schedule outcome; no consumption in the case of a too high electricity price, and full load consumption in the case of a too low electricity price. The prices that fall within the green and gray value range, however, might lead to different schedules depending on the submitted bid curve. We therefore mark this range in shaded red. More than 60\% of the hours have a day-ahead price within this interest range. 

The RES data is based on the day-ahead forecast and realized production of the total onshore wind production in the DK1 zone, taken from the public platform \textit{Energi Data Service}~\citep{EnerginetForecastResolution}. While the electrolyzer would likely time-match with a specific  individual renewable farm or farms, we have chosen to use this aggregated data as it is publicly available. Aggregated RES production is typically less volatile than on individual plant level due to spatial smoothing effects, and is likely to have smaller forecast errors than individual farms \citep{Focken2002Short-termEffects}. Hence, the effects of uncertainty as seen in this work is more likely underestimated than overestimated due to the use of the aggregated wind data. Figure \ref{fig:wind_distribution} shows the forecast and associated forecast errors on the year 2024, given in per unit of maximum production. The wind forecast error is computed as
\begin{align}
    \epsilon_t =  \tilde{P}^{\rm w}_t - \hat{P}^{\rm w}_t,
\end{align}

where $\tilde{P}^{\rm w}_t$ is the realized wind production and $\hat{P}^{\rm w}_t$ the forecasted production of an hour. From Figure \ref{fig:wind_distribution} we see that there is a tendency of positive forecast error, i.e., higher realization than forecast, at the mid-to-high forecast values. The available data includes only a point-forecast and a realization. Information such as the confidence interval of the prediction is not available in the obtained dataset, and is typically only provided from merchant forecast providers. 
In the following results, we sample scenarios based on the realized wind production in previous times with a forecasted value similar to the current.
We sample nine scenarios and include the day-ahead point forecast as the last scenario, for a total of $N = 10$ scenarios. We refer to this method as the \textit{K-nearest neighbors} (KNN) method. The sampling method is described in detail in Appendix~\ref{appenidx:sampling_and_parameters}. 

The base case inputs are specified in Table \ref{tab:case_study_base_inputs} of Appendix~\ref{appenidx:sampling_and_parameters}.  We perform a sensitivity analysis on the number of scenarios in the later Section \ref{sec:sensitivity_sampling}, and on the gray and green hydrogen price and subsidy in Section \ref{sec:h2_price_sensitivity}.

\subsection{Deriving electrolyzer bid curves}

For every hour of the year, we compare \textbf{three general bid curves:} \textit{i}) a bid curve based on the day-ahead point forecast of the wind production, \textit{ii}) a bid curve based on scenarios of wind-production derived from the day-ahead forecast and, \textit{iii}) a bid curve based on perfect information of the realized wind production.  In the following, we summarize the derivation of the bid curves. A detailed derivation of these curve is provided in Appendix \ref{bid_curve_derivation}. 

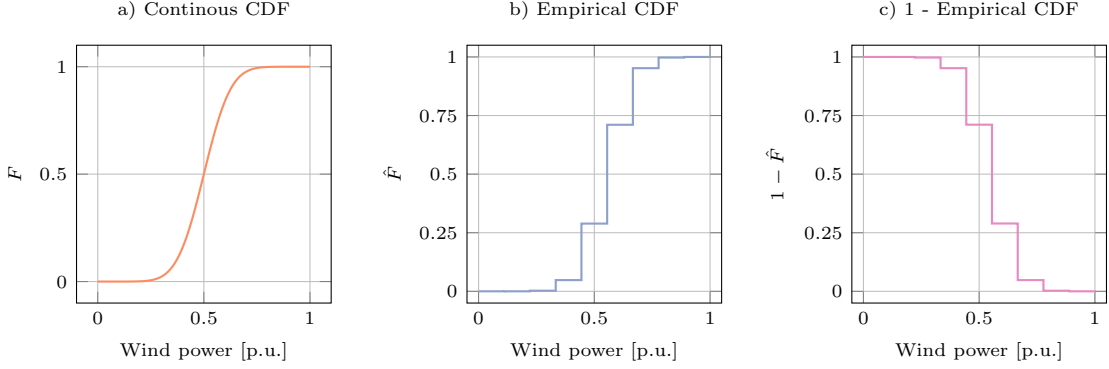
\begin{figure}[t]
\centering

\begin{subfigure}[t]{0.32\textwidth}
\centering
\begin{tikzpicture}
\begin{axis}[
    width=\linewidth,
    height=5cm,
    xlabel={Wind power [p.u.]},
    ylabel={$F$},
    domain=0:1,
    samples=100,
    grid=both,
    title={\scriptsize a) Continous CDF },
    tick label style={font=\scriptsize},
    label style={font=\scriptsize},
]
\addplot[color2, thick] {gausscdf(0.5, 0.1)};
\end{axis}
\end{tikzpicture}
\end{subfigure}
%
\begin{subfigure}[t]{0.32\textwidth}
\centering
\begin{tikzpicture}
\begin{axis}[
    width=\linewidth,
    height=5cm,
    xlabel={Wind power [p.u.]},
    ylabel={$\hat F$},
    xmin=-0.05, xmax=1.05,
    ymin=-0.05, ymax=1.05,
    samples=10,
    domain=0:1,
    ytick={0,0.25,0.5,0.75,1},
    xtick={0,0.5,1},
    grid=both,
    title={\scriptsize b) Empirical CDF},
    tick label style={font=\scriptsize},
    label style={font=\scriptsize},
]
\addplot+[const plot mark left, no markers, color3, thick] 
    {gausscdf(0.5, 0.1)};
\end{axis}
\end{tikzpicture}
\end{subfigure}
%
\begin{subfigure}[t]{0.32\textwidth}
\centering
\begin{tikzpicture}
\begin{axis}[
    width=\linewidth,
    height=5cm,
    xlabel={Wind power [p.u.]},
    ylabel={$1 - \hat F$},
    xmin=-0.05, xmax=1.05,
    ymin=-0.05, ymax=1.05,
    samples=10,
    domain=0:1,
    ytick={0,0.25,0.5,0.75,1},
    xtick={0,0.5,1},
    grid=both,
    title={\scriptsize c) 1 - Empirical CDF},
    tick label style={font=\scriptsize},
    label style={font=\scriptsize},
]
\addplot+[const plot mark left, no markers, color4, thick] 
    {1 - gausscdf(0.5, 0.1)};
\end{axis}
\end{tikzpicture}
\end{subfigure}
\caption{Comparison between the continuous CDF, the empirical CDF and 1 - the empirical CDF for wind power production assuming a normal distribution with mean 0.5 and standard deviation 0.1.}
\label{CDFs}
\end{figure}

The bid curve should consist of price-quantity pairs, where the bid price is given by the expected marginal value of the hydrogen production at this quantity.  In the continuous case for any bid quantity $q$, we have:
\begin{align}
    {\rm price}(q) = \eta \Big (\pi^{\rm gray} + \pi^{\rm green} \big(1-{\rm F}_{P^{\rm w}}(q) \big)\Big), \qquad P^{\rm w}: \Omega,
\end{align}
where ${\rm F}( q)$ is the cumulative distribution of the wind production, as illustrated in Figure \ref{CDFs} a).  The factor $\big(1-{\rm F}_{P^{\rm w}}(q) \big)$ is interpreted as the probability of the wind realization being higher than a certain quantity $q$ for a given distribution of the wind production $\Omega$. We multiply this factor with the green subsidy $\pi^{\rm green}$, as the subsidy is only awarded for electrolyzer consumption up to the RES production. Further, all consumption yields the the gray hydrogen price $\pi^{\rm gray}$, and all consumption is multiplied with the electrolyzer efficiency $\eta$. 

To bid into the electricity market, the bid curve should be discretized. We therefore draw a set of scenarios $\mathcal{S}$ for wind production from the continuous distribution $\Omega$, and evaluate the empirical cumulative distribution function ${\rm \hat{F}}(q)$:
\begin{align}
    {\rm price}(q) = \eta \Big (\pi^{\rm gray} +\pi^{\rm green} \big ( 1 - {\rm \hat{F}}_S(q) \big ) \Big ), \qquad s \in \mathcal{S}.
    \label{eq:emirical_step_curve}
\end{align}
The empirical cumulative distribution function ${\rm \hat{F}}_S(q)$ is interpreted as the summed probability of scenarios that fall below the quantity $q$, and is illustrated in Figure \ref{CDFs} b). Hence, $(1- {\rm \hat{F}}_S(q))$ represents summed probability of scenarios that are higher than $q$, as illustrated by Figure \ref{CDFs} c).
This results in discrete and decreasing price-quantity steps, where a new (and lower) step is taken for every scenario of wind production $P^{\rm w}_s$.

\begin{figure}[ht]
\centering
\begin{minipage}[t]{0.33\textwidth}
\flushleft
\scriptsize
\begin{tikzpicture}
\begin{axis}[
    width=\linewidth,
    height=5cm,
    title={\textit{i}) Point forecast bid},
    title style={font=\scriptsize},
    xlabel={Quantity [MWh]},
    xlabel style={font=\scriptsize},
    ylabel={Price  [EUR/MWh]},
    ylabel style={font=\scriptsize},
    tick label style={font=\scriptsize},
    grid=both,
    grid style={dotted},
    legend style={font=\scriptsize, at={(0.5,-0.3)}, anchor=north, legend columns=-1},
    ymin=0, ymax=80,
    xmin=0, xmax=110,
    xtick={0, 25, 50, 75, 100},
    xticklabels = {0, , $\hat P^{\rm W}$, , $P^{\rm h}_{}$},
    ytick={20, 60},
    yticklabels={$ \lambda^{\rm gray} $, $ \lambda^{\rm green} $}
]
\addplot[color=color2, thick] coordinates {(0,60) (50,60)};
\addplot[color=color2, thick] coordinates {(50,60) (50,20)};
\addplot[color=color2, thick] coordinates {(50,20) (100,20)};
\addplot[color=color2, thick] coordinates {(100,20) (100,0)};
\end{axis}
\end{tikzpicture}
\end{minipage}
\begin{minipage}[t]{0.32\textwidth}
\centering
\scriptsize
\begin{tikzpicture}
\begin{axis}[
    width=\linewidth,
    height=5cm,
    title={\textit{ii}) Scenario forecast bid},
    title style={font=\scriptsize},
    xlabel={Quantity [MWh]},
    xlabel style={font=\scriptsize},
    ylabel={},
    ylabel style={font=\scriptsize},
    tick label style={font=\scriptsize},
    grid=both,
    grid style={dotted},
    legend style={font=\scriptsize, at={(0.5,-0.3)}, anchor=north, legend columns=-1},
    ymin=0, ymax=80,
    xmin=0, xmax=110,
    xtick={0, 25, 50, 75, 100},
    xticklabels = {0, $\hat P^{\rm W}_{1}$, , $\hat P^{\rm W}_{2}$, $P^{\rm h}_{}$},
    ytick={20, 60},
    yticklabels={}
]
\addplot[color=color3, thick] coordinates {(0,60) (25,60)};
\addplot[color=color3, thick] coordinates {(25,60) (25,40)};
\addplot[color=color3, thick] coordinates {(25,40) (75,40)};
\addplot[color=color3, thick] coordinates {(75,40) (75,20)};
\addplot[color=color3, thick] coordinates {(75,20) (100,20)};
\addplot[color=color3, thick] coordinates {(100,20) (100,0)};
\end{axis}
\end{tikzpicture}
\end{minipage}
\begin{minipage}[t]{0.32\textwidth}
\centering
\scriptsize
\begin{tikzpicture}
\begin{axis}[
    width=\linewidth,
    height=5cm,
    title={\textit{iii}) Perfect information bid},
    title style={font=\scriptsize},
    xlabel={Quantity [MWh]},
    xlabel style={font=\scriptsize},
    tick label style={font=\scriptsize},
    grid=both,
    grid style={dotted},
    legend style={font=\scriptsize, at={(0.5,-0.3)}, anchor=north, legend columns=-1},
    ymin=0, ymax=80,
    xmin=0, xmax=110,
    xtick={0,  50, 70,  100},
    xticklabels = {0,  $\hat P^{\rm W}$, $\tilde P^{\rm W}$,  $P^{\rm h}_{}$},
    ytick={20, 60},
    yticklabels={}
]
\addplot[color=color1, thick] coordinates {(0,60) (70,60)};
\addplot[color=color1, thick] coordinates {(70,60) (70,20)};
\addplot[color=color1, thick] coordinates {(70,20) (100,20)};
\addplot[color=color1, thick] coordinates {(100,20) (100,0)};
\end{axis}
\end{tikzpicture}
\end{minipage}

\caption{Illustration of the three general electrolyzer day-ahead bid curves.}
\label{fig:illustration_bid_curves}
\end{figure}
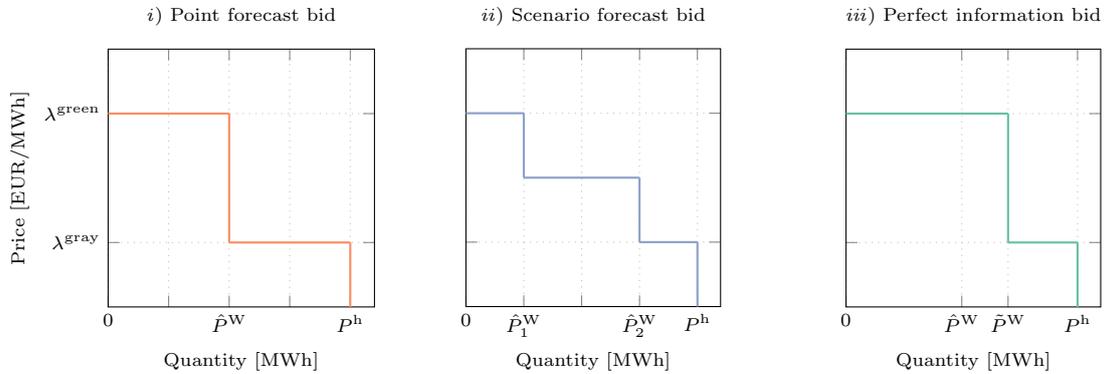

Figure \ref{fig:illustration_bid_curves} illustrates the three general bid curves. The height of a step, i.e., the bid price, is given by eq. (\ref{eq:emirical_step_curve}). 
The point forecast curve, Figure \ref{fig:illustration_bid_curves} \textit{i}), and the perfect information curve, Figure \ref{fig:illustration_bid_curves} \textit{iii}), are both characterized by two steps. The higher step represents the green hydrogen value $\lambda^{\rm green}$, i.e., when a subsidy is received as given by eq. (\ref{green_price}). The length of the step is given by the amount of wind production, either forecasted $\hat P^{\rm W}$,  or realized $\tilde P^{\rm W}$.  The lower step represents the gray hydrogen value, i.e., the remaining consumption up the the electrolyzer capacity $ P^{\rm h}$, which does not qualify for a subsidy. For both these curves, this binary nature of the curve arises from the fact that you are either 100\% certain that you will receive the subsidy, or 100\% certain that you will not. This is equivalent to drawing a single scenario in eq. (\ref{eq:emirical_step_curve}), equal either to the point forecast or the realized production. Note that in a case where the wind forecast or realization is larger than the electrolyzer capacity $P^{\rm h}$, only the first step will appear on the bid curve, and vice versa in case of a zero forecasted or realized production.

The scenario forecast curve, Figure \ref{fig:illustration_bid_curves} \textit{ii}), however, can have two or more steps depending on the number of wind production scenarios included in its derivation, with the number of steps always equal to $S+1$, where $S$ is the number of scenarios included with a forecast lower than the electrolyzer capacity.  The height of the steps, i.e., the price bids, is a function of the probability of the wind production being at least the corresponding quantity. In the example curve, we have two scenarios with equal probability. Hence, there is a 100\% chance that the wind production will be equal to or greater than the lower value scenario. Further, the probability of the wind production being between the lower and higher scenario is 50\%, which yields a bid price equal to receiving half of the green compliance subsidy. Finally, the probability of the wind production being larger than the highest scenario is 0\%, hence the bid price is equal to the gray hydrogen value.

While the three methods of deriving the three bid curves are constant, the input information regarding the wind production changes per hour, which means that the derived bid curves change per hour.

\subsection{Electrolyzer profits comparison}
\begin{figure}[t!]
\begin{minipage}[t]{0.45\textwidth}
    \centering
    \includegraphics[width=\linewidth]{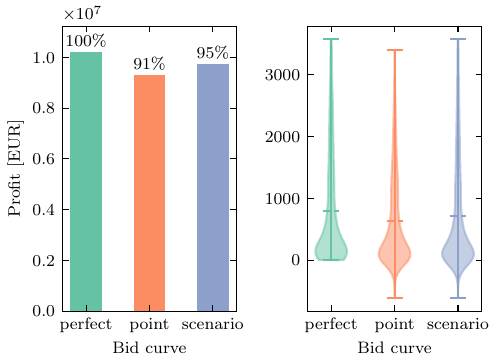}
    \caption{Total (left) and distribution of (right) profits over 2024 of the perfect information, scenario based and point forecast curves.}
    \label{fig:Profit_total}
\end{minipage}
\hfill
\begin{minipage}[t]{0.45\textwidth}
    \centering
    \includegraphics[width=\linewidth]{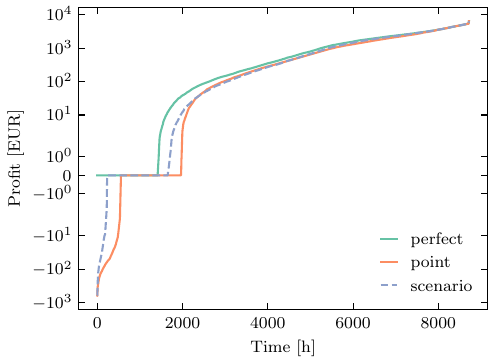}
    \caption{Profit duration curve with logarithmic y-axis.}
    \label{fig:Profit_duration}
\end{minipage}
\end{figure}

We first compare the three bid curves in terms of operational profits. Figure \ref{fig:Profit_total} shows the total profits of each curve relative to the perfect information curve (left) as well as the distribution of profits over the year (right) for each curve. The point forecast curve achieves 91\% of the perfect information curve profit in total over the year 2024. This is increased to 95\% for the scenario based bid curve. The distribution of profits, Figure \ref{fig:Profit_total} (right) shows that while the perfect information curve achieves a non-negative profit in all hours, both the point forecast curve and scenario based curve have instances of negative profits. 

To better illustrate the distributional differences between especially the point forecast and scenario based curves, we plot their respective profit duration curves on a logarithmic axis in Figure \ref{fig:Profit_duration}, where the profit of each hour is sorted in ascending order. The figure further illustrates that both the point forecast and scenario based curves achieve negative profits in some hours. However, the point forecast curve have about twice the amount of negative profit hours compared to the scenario based curve. This is where the scenario based curve gains the most advantage compared to the point forecast curve.


\begin{figure}[t]
    \centering
    \begin{subfigure}[t]{0.48\textwidth}
        \centering
        \includegraphics[width=\linewidth]{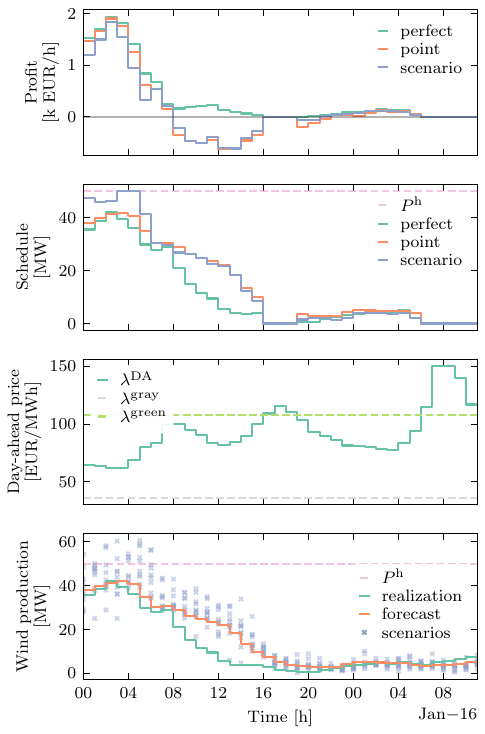}
        \caption{Example day where neither of the forecasted curves results in a schedule that matches with the realized wind production. The wind realization is lower than the forecast.}
        \label{fig:profits_per_hour_1}
    \end{subfigure}%
    \hfill
    \begin{subfigure}[t]{0.48\textwidth}
        \centering
        \includegraphics[width=\linewidth]{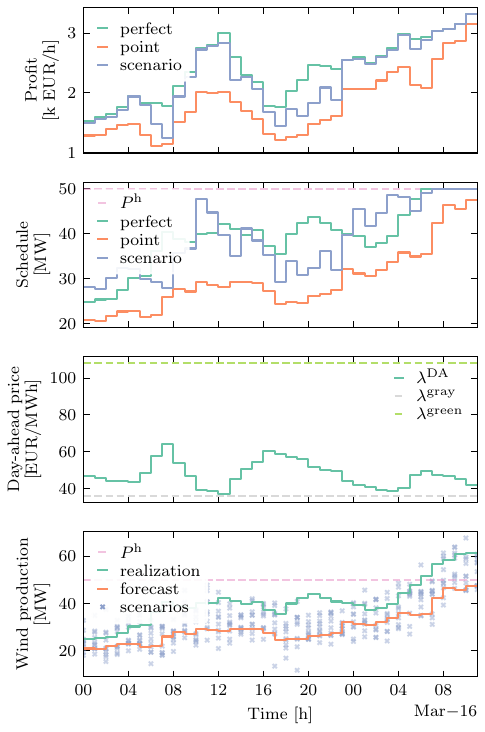} 
        \caption{Example day where the scenario curve consistently out-performs the point forecast.}
        \label{fig:profits_per_hour_2}
    \end{subfigure}
    \caption{Profit and schedule of each curve, day-ahead price and wind production (realized and forecasted) over 36 hours. $\lambda^{\rm DA}$ notes the day-ahead price, $\lambda^{\rm gray}$ the gray hydrogen value, $\lambda^{\rm green}$ the green hydrogen value, and $P^{\rm h}$ the electrolyzer capacity.}
    \label{fig:profits_comparison}
\end{figure}

We next investigate the profits of the three curves during specific hours. Figure \ref{fig:profits_comparison} shows the resulting electrolyzer profits (top), electrolyzer schedule (second) , the day-ahead price (third) and wind production data (bottom) for two selections of 36 hour periods in 2024. Figure \ref{fig:profits_per_hour_1} shows an example where the neither the point forecast nor scenarios capture the realized wind production between the hours 8-16 well. This leads to negative profit in those hours for both curves, with a slightly lower loss for the scenario based curve. On the contrary, Figure \ref{fig:profits_per_hour_2} shows an example where the scenario based curve consistently outperforms the point forecast curve, with instances of achieving the near perfect profit. In this example, the realized wind production was consistently higher than the forecast, while the day-ahead prices where in the lower range, but not below the gray hydrogen price. The relatively low day-ahead prices means that the scenario based curve yields a schedule which is higher than the point forecast, and is therefore better aligned with the perfect information schedule.



\subsection{Following the wind production}
The goal of the green subsidy is to incentivize the electrolyzer to match its consumption to an (additional and spatially correlated) RES. We therefore compare the electrolyzers scheduled consumption against the realization of the wind production. 
We refer to the case where the electrolyzer consumes more than the realized as ``gray" consumption, $P^{\rm gray}$, which we compute as follows:
\begin{align}
    P^{\rm gray}_t = \text{max}(q^{\rm curve*}_t - \tilde{P}^{\rm w}_t, 0) .
\end{align}
Here, $q^{\rm curve*}_t$ is the cleared consumption quantity of the electrolyzer in an hour for a given bid curve.
As expected, the perfect information curve best follows the wind production, as it would only lead to a gray consumption in hours where the day-ahead price is equal to or less than the gray hydrogen value. The point forecast and scenario forecast curves might, however, lead to a gray consumption due to the realized wind production being lower than expected. 

\begin{figure}
    \centering
    \includegraphics[width=0.5\linewidth]{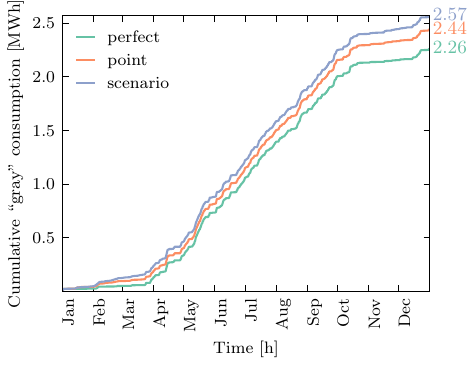}
    \caption{Cumulative sum of the electrolyzers over-consumption compared to the realized wind production, i.e., ``gray" consumption.}
    \label{fig:wind_deficit}
\end{figure}

While we have seen that the scenario forecast curve improves the electrolyzers profits, Figure \ref{fig:wind_deficit} indicates that it does not improve the electrolyzer's ability to follow the realized wind production. The point forecast curve leads to a 8\% increase in gray consumption compared to the perfect information curve, while the scenario forecast curve leads a 25\% increase. These increases cannot be attributed to consumption at low day-ahead prices, i.e., during hours when producing gray hydrogen is profitable, as they are captured by the perfect information curve. Instead, these cases of  gray consumption are due to the realized wind production being lower than the what the point forecast or scenario forecast indicates. 

Table \ref{tab:hydrogen_production} shows the total amounts of green and gray hydrogen produced over the year for each of the three curves. Here, green hydrogen is that produced from consumption that time-matched with the realized wind production, while the gray hydrogen is that of the remainder. The share of green and gray hydrogen compared to the total amount produced for each curve is noted. This share remains similar under each case, with the perfect forecast curve share slightly higher than the point forecast and scenario based curves. 

\begin{table}[ht]
\small
\centering
\caption{Total green and gray hydrogen production.}
\begin{tabular}{lccccc}
\hline
Method & Green [t] & Green [\%] & Gray [t] & Gray [\%] & Total [t] \\
\hline
Perfect  & 2714.02 & 78.8 & 731.50 & 21.2 & 3445.52 \\
Point    & 2396.59 & 75.2 & 790.16 & 24.8 & 3186.75 \\
Scenario & 2504.65 & 74.5 & 886.65 & 25.5 & 3391.30 \\
\hline
\end{tabular}
\label{tab:hydrogen_production}
\end{table}



This result highlights a potential pitfall of time-matching regulation due to the uncertain nature of RES. While the point forecast curve also leads to an increase in wind deficit compared to the electrolyzer consumption, this is intensified when the electrolyzer operates in an uncertainty aware manner. The further impact of this result on the system is not easily quantified within the framework of the current study. 

Our analysis of the electrolyzer's profit indicates that the scenario forecast curve outperforms the point forecast curve on realizations that lead to negative profits. These cases are characterized by day-ahead prices in the higher end of the relevant price range, and large wind forecast errors. High day-ahead prices could be an indication of scarcity of cheaper generation resources such as RES. Hence, it might be an advantage that the scenario based curve performs better, and more closely follows the realized wind production, during these hours in particular.

\subsection{Sensitivity analysis}

\subsubsection{Sampling}
\label{sec:sensitivity_sampling}
Next, we perform a sensitivity analysis of the resulting profits using various numbers of samples. The scenario generation method applied is described in Appendix~\ref{appenidx:sampling_and_parameters}. 

\begin{figure}[h]
    \centering
    \includegraphics[width = \textwidth]{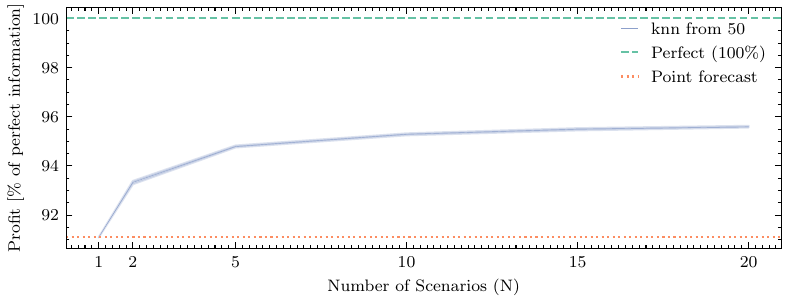}
    \caption{Profits as a percentage of the perfect information case with variations of sampling approach. }
    \label{fig:profits_sampling_var}
\end{figure}

Figure~\ref{fig:profits_sampling_var} shows the resulting profits as a percentage of the perfect information case, which is indicated at 100\% by the blue dashed line. The point forecast profit is indicated by the orange dotted line. All sampling methods include the current point forecast as a scenario. Hence, for a single sample ($N = 1$), the profit is always equal to that of the point forecast curve. The results show a diminishing increase in performance with increasing number of samples, flattening out at 10 samples. The highest profit is achieved with 20 samples, with 95.6\% of the perfect information profit. 

Figure~\ref{fig:profits_sampling_var} further shows the results cross-validation on the sample selection as a band around the mean line. For each number of samples $N$, we have drawn 10 different sets of scenarios following the methodology described in Appendix~\ref{appenidx:sampling_and_parameters}. However, the variation is so small, that the band nearly blends with the mean line. In conclusion, the increase in profits achieved by the scenario based bid curve compared to the point forecast curve appears to be stable.



\subsubsection{Hydrogen prices and subsidies}
\label{sec:h2_price_sensitivity}
\begin{figure}[h]
    \centering
    \includegraphics[width = \textwidth]{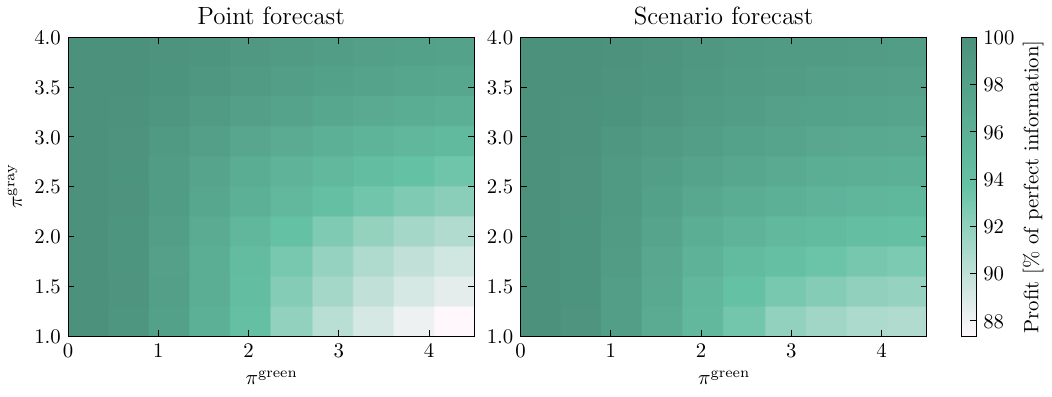}
    \caption{Profit sensitivity towards the gray and green hydrogen prices.}
    \label{fig:profit_sensitivity}
\end{figure}

Further, we conduct a sensitivity analysis on the gray hydrogen price $\pi^{\rm gray}$ and the green hydrogen subsidy $\pi^{\rm green}$. Figure \ref{fig:profit_sensitivity} shows the sensitivity on the profits of the point forecast and scenario based methods, relative to the perfect information profit. Both curves perform relatively better under high gray hydrogen prices and low green hydrogen subsidies. Additionally, the difference in profits between the two methods decrease. This is expected, as with a lower green hydrogen subsidy, the effect of having imperfect information decreases. Further, Figure \ref{fig:profit_sensitivity} shows that if the gray hydrogen price increases, the effect of imperfect information decreases. This is due to two reasons. First, with higher gray hydrogen prices, the relative gain of the green hydrogen subsidy decreases. Second, the higher gray hydrogen price will move the interest range of day-ahead prices, as defined in Figure \ref{fig:price_duration}, further up on the y-axis. As there are less hours that fall within this higher interest range, there are less hours with the possibility of different scheduling decision depending on the bid curve method. However, as commercial gray hydrogen is typically a derivative of natural gas, its production cost is highly dependent on the natural gas price \citep{InternationalEnergyAgency2024Global2024}. While there are still gas fired generation in the continental European power system, it is plausible that an increase in the natural gas price would yield an increase in electricity prices. If so, that could mean larger effects of the RES uncertainty also at higher gray hydrogen prices.

Next, we investigate the sensitivity towards the gray hydrogen price and green subsidy on the share of green hydrogen produced by the electrolyzer, i.e., the hydrogen that is produced with consumption that time-matches with the realized wind production. For both the point forecast and the scenario forecast bid curves, the share of green hydrogen increases with decreasing gray hydrogen price and increasing green hydrogen subsidy. When the gray hydrogen price decreases, the number of hours where the electrolyzer is willing to produce without gaining the time-matching subsidy decreases. Thus, with higher hydrogen subsidy, the electrolyzer will produce more during hours where the green hydrogen subsidy is available, and the share of green hydrogen production increases.  

\begin{figure}[h]
    \centering
    \includegraphics[width = \textwidth]{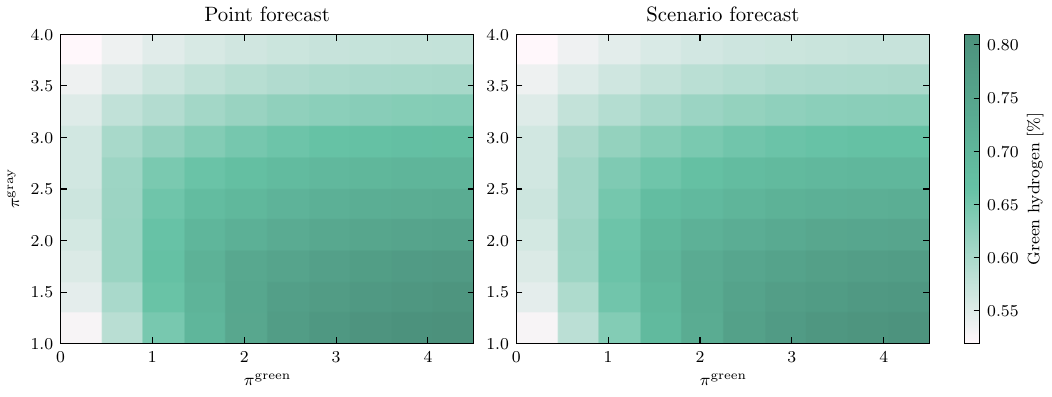}
    \caption{Green hydrogen sensitivity towards the gray and green hydrogen prices}
    \label{fig:green_sensitivity}
\end{figure}

\subsubsection{The electrolyzer's impact on clearing prices}

A inherent problem with using a historically cleared price to compute the profits of the electrolyzer is that this approach cannot capture the electrolyzer's part in the price formation process. When the electrolyzer participates in the market, the demand will generally increase, which might increase the clearing prices. This could in turn reduce the electrolyzer's profits.

In order to evaluate the potential impact of the electrolyzer on the historical prices, we employ historical aggregate demand and supply curves from 2024. These curves are annonymized representations of the real, historical supply and demand curves, aggregated on the national Danish level~\citep{NordPoolAggregatedDay-ahead}. 

To compute the electrolyzer's impact on the cleared price, we first find the cleared prices without any electrolyzer present. This price is given by the merit-order of the historical aggregated curves, i.e., the price at the point where the two curves intercept. Then, the bid curve of the electrolyzer is computed for every hour, and placed in the merit order of the aggregate demand curve of the that hour. The new clearing price is found where the updated demand curve and supply curve intercepts. 

\begin{figure}[t]
    \centering
\includegraphics[]{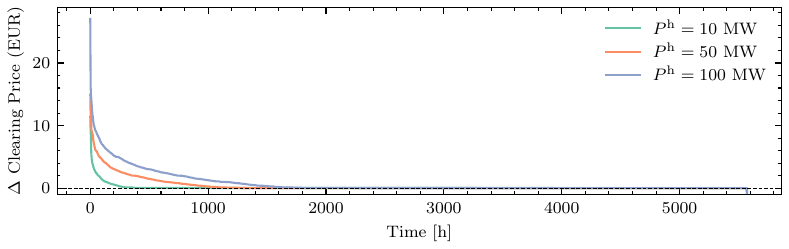}
    \caption{The impact of various electrolyzer capacities on the cleared electricity price in Denmark, bidding with the perfect information curve.}
    \label{fig:price_impact}
\end{figure}

Figure~\ref{fig:price_impact} shows the resulting changes in cleared electricity price for three electrolyzers with capacities of 10, 50 and 100~MW respectively. All electrolyzers are modeled with perfect information bid curves. The change in clearing price is sorted in descending order in each case, and reveals that an electrolyzer capacity of 10~MW only changes the price in about 250 hours, while the 100~MW electrolyzer has an impact in about 1600 hours. While the full year consist of 8760 hours, only 5400 hours are included in the plot, due to missing data. While the larger electrolyzer capacities impact the cleared prices in a significant amount of hours, the magnitude of this change stays below EUR 5 per MWh in most hours. 

Note that while these results might give an indication of the electrolyzer's impact on the clearing price, they do not capture, for instance, the effect of congestion in the transmission system which could further impact the changed is cleared prices. 

\section{Conclusions and future work}
In this work we derive a bid curve for an electrolyzer to participate in electricity markets when operating under green hydrogen regulation. Specifically, we address the temporal-matching requirement. Under this requirement, the uncertainty of the matching RES propagates to the value of hydrogen for the electrolyzer. Therefore, we develop a bid curve for the electrolyzer that takes this renewable production uncertainty into account. 

We formulate a linear program which maximizes the expected profit of the electrolyzer given an electricity price and scenarios of the matching renewable production. Further, we derive the bid curve from the KKT-conditions of this linear program. The resulting curve is a stepwise, price-quantity curve, where the height of a step represents the probability of receiving the subsidy for the  corresponding quantity.

We show that bidding with an uncertainty-aware curve outperforms bidding with a point-forecast curve in terms of profits, even under naive assumptions about the renewable production scenarios. Where the point forecast bid curve achieves 91\% of the perfect information profit, the scenario based curve yields 95\%. However, the uncertainty-aware bid curve also leads to an increase in electrolyzer consumption that exceeds the production of the matching RES, compared to the point forecast curve. While the current study cannot conclude on what effects this increase might have on the power system emissions, it points to the fact that uncertainty might distort some of the intended effects of the green hydrogen subsidy and temporal regulation. 

This working paper leaves two important directions for future work. The analysis conducted on the historical price signals could be extended to the historical aggregate supply and demand curves. This could provide insights into how the electrolyzer's impact on the cleared price might change the relative performance of the different curves. For example, when using the historical aggregate curves, it is no longer guaranteed that the perfect information bid curve will lead to the highest profits, as it might result in a cleared electricity price higher than in other cases. 

The impact of uncertainty on the effectiveness of temporal regulations should be investigated. The bid curves derived in the current work could be used to represent the electrolyzer in long-term studies, in order to fully capture the impacts of uncertainty of the efficacy of green hydrogen subsidies and temporal regulation.


\newpage

\bibliographystyle{elsarticle-harv} 
\bibliography{ref_test}
\newpage
\appendix
\section{Scenario generation method and case study parameters}
\label{appenidx:sampling_and_parameters}

\begin{table}[ht]
\centering
\caption{Case study base parameters.}
\label{tab:case_study_base_inputs}
\begin{tabular}{lll} \hline
Gray hydrogen price $\pi^{\rm gray}$              & 2                     & EUR/kg               \\
Green hydrogen subsidy $\pi^{\rm green}$            & 4                     & EUR/kg               \\
Electrolyzer efficiency  $\eta$         & 18                    & kg/MWh             \\
Electrolyzer capacity  $P^{\rm h}$           & 50                    & MW \\ \hline
Electricity price data            & Day-ahead DK1 2024    & EUR/MWh              \\
RES data                          & Onshore wind DK1 2024 & MWh                \\
RES capacity                      & 66                  & MW            \\
RES scenarios $N$                      & 10                  & -                \\
Error sample pool $K$                      & 50                  & -                \\\hline
\end{tabular}
\\ 
\end{table}

Table~\ref{tab:case_study_base_inputs} provides the input parameters used in the case study presented in this work. The electrolyzer capacity is set to 50~MW, and we assume that this capacity is 75\% of that of the matching wind farm. The wind farm therefore has a capacity of 66~MW. The optimal sizing between the electrolyzer and wind farm capacities has not been addressed. The sampling method used to create scenarios of wind production is described in the following.

The available data provides a forecast and realization of the total onshore wind production in DK1 \citep{EnerginetForecastResolution}. The data is scaled to the wind farm capacity of 66~MW. To create scenarios of potential forecast errors, we utilize historical data from the years 2023 and 2024. For any hour $t$, we find the $K = 50$ hours with forecasted value closest to that of the current. Hence, we achieve the set \textit{K-nearest neighbors} to the current forecast value. For this set of 50 hours, we compute the per unit forecast error. This gives us a set of 50 per unit errors that we multiply with the current forecast. From the resulting set of 50 errors, we randomly draw a number of $N-1$ scenarios, where in the base case $N=10$. The $N$\textsuperscript{th} scenario included is the current forecast value, i.e., assuming a forecast error equal to zero. Recall that the case study only considers year 2024, hence the 2023 data is only applied to create scenarios. Only data that precedes the current hour is used, and a real electrolyzer operator would therefore have access to this data. 

The scenario generation method is mathematically described as follows. Let \(\hat P^{\rm RES}_t\) denote the forecasted renewable output at hour \(t\).  Define the historical index set,
\[
  \mathcal{T}_t \;:=\;\{\,j : j < t, \quad t \in \mathcal{T}^{\mathrm{full}}\},
\]
to include only past hours, where $\mathcal{T}^{\mathrm{full}}$ is the full set of hours in the dataset.  From \(\mathcal{T}_t\), select the $K$ indices whose forecasts are closest to \(\hat P^{\rm RES}_t\):
\[
  \mathcal{J}_t \;:=\;{\arg\min}_{\substack{\mathcal{J}\subset\mathcal{T}_t,\,|\mathcal{J}|=K}}
    \sum_{j\in\mathcal{J}} \bigl|\hat P^{\rm RES}_j - \hat P^{\rm RES}_t\bigr|.
\]

The resulting set can also be described as:
\[
  \mathcal{J}_t \;:=\;\bigl\{\,j : \hat P^{\rm RES}_j \text{ is among the $K$ past closest values to } \hat P^{\rm RES}_t\,\bigr\}.
\]

 In the base case, $K=50$.
For each \(j \in \mathcal{J}_t\), compute the per‑unit forecast error,
\[
  e_j \;=\;\frac{\tilde P^{\rm RES}_j - \hat P^{\rm RES}_j}{\hat P^{\rm RES}_j}\,,
\]
where \(\tilde P^{\rm RES}_j\) is the realized renewable output at hour \(j\).  This yields a set of errors \(\{e_j : j \in \mathcal{J}_t\}\).  

Form the scenario set by uniformly sampling \(N\) errors (without replacement) from \(\{e_j : j \in \mathcal{J}_t\}\).  Denote the sampled errors by \(\{e_s\}_{s=1}^{N-1}\).  Then define the scenario values of the green signal, $P_s^{\rm RES}$, as
\begin{align*}
  & P_s^{\rm RES}
  = \hat P^{\rm RES}_t \,\bigl(1 + e_s\bigr)
  \quad & \text{for }s = 1,\dots,N-1, \\ 
  &P_s^{\rm RES}
  = \hat P^{\rm RES}_t 
  \quad &\text{for }s = N, \\
\end{align*}
where each scenario is assigned an equal probability \(\rho_s = 1/N\).

\newpage

\section{Scenario-based price-quantity bid curve derivation}
\label{bid_curve_derivation}

To represent the uncertainty on the green signal $P^{\rm RES}$, we consider a discrete set of scenarios $\{P_s^{\rm RES}\}_{s \in \mathcal{S}=\{1,..,S\}}$, associated with the probabilities $\rho_s = \mathbb{P}(P^{\rm RES} = P^{\rm RES}_{s}) \in ]0,1]$ where $\sum_{s \in \mathcal{S}} \rho_s = 1$. Without any loss of generality, we assume that these scenarios are unique and ordered, such that $0 \leq P^{\rm RES}_1 < P^{\rm RES}_2 < ... < P^{\rm RES}_{s-1} < P^{\rm RES}_{s} < ... < P^{\rm RES}_S \leq P^{\rm h}$.\footnote{Note that, in the case where the green signal $P^{\rm RES}$ represents a value that can exceed the maximum electrolyzer load, we can transform into a signal that is saturated at $P^{\rm h}$.}

The electrolyzer's optimal stepwise bid curve in the DA market can be represented by a discrete number $K$ of non-decreasing price-quantity pairs $b^{\rm DA}_k = \left(\lambda^{DA}_k,P^{DA}_k\right)$ for $k =1,...,K$. In the following we detail the methodology to compute the pieces of this optimal bid curve, given a discrete set of green signal scenarios.

\subsection{Expected value of hydrogen}

\begin{figure}[!ht]
\centering
\label{perfect_forecast}
\begin{tikzpicture}
    \begin{axis}[
        width=6cm,
        height=6cm,
        xlabel={Energy [MWh]},
        ylabel={Hydrogen value [\$]},
        grid=both,
        grid style={dotted},
        legend style={at={(0.5,-0.3)},anchor=north,legend columns=-1},
        ymin=0, ymax=200,
        xmin=0, xmax=110,
        xtick={0, 25, 50, 75, 100},
        ytick={0, 25, 50, 75, 100, 125, 150, 175 },
        yticklabels = {0, },
        xticklabels = {0, $P^{\rm RES}$, , , $P^{\rm h}$}
    ]
        \addplot[ color=green, thick] coordinates {(0,0) (25,25*2)};
        \addlegendentry{${\rm g} (x)$};
        \addplot[forget plot, dashed, color=green, thick] coordinates {(0,0) (100,100*2)};
        \addplot[color=gray, opacity = 0.7, thick] coordinates {(25,25*2) (100,75*1+25*2)};
        \addlegendentry{${\rm f} (x)$};
        \addplot[forget plot, dashed, color=gray, opacity = 0.7, thick] coordinates {(0,25) (25,25*2)};
        ;
    \end{axis}
\end{tikzpicture}
\caption{The piece-wise linear hydrogen value curve. The binding pieces are marked by full lines, and the dashed lines mark the non-binding parts of the functions, illustrating that they form a concave hull.}
\label{fig:concave_hull_h2_value}
\end{figure}
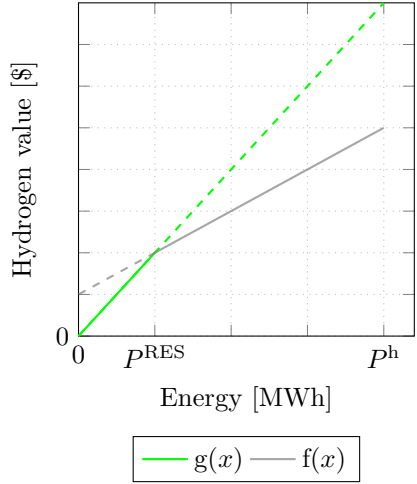

Let us first define the value of hydrogen for a given value of the green signal, $P^{\rm RES}$, as a function of the day-ahead purchase quantity. The hydrogen value function is a two-piece linear curve, with the first piece representing the value under the green premium, and the second piece representing the value of gray hydrogen when exceeding the green-signal, as illustrated by Figure \ref{fig:concave_hull_h2_value}.

Let's call the steeper, green piece on the curve ${\rm g} (x)$ and the flatter, gray piece on the curve ${\rm f} (x)$. We define the value of hydrogen as the following function:
\begin{align}
    {\rm \gamma} : x \rightarrow \begin{cases}
        {\rm g} (x) \quad \text{if} \quad x \leq P^{\rm RES}, \\
        {\rm f} (x) \quad \text{if} \quad P^{\rm RES} < x \leq P^{\rm h}.
    \end{cases}
\end{align}
As the gray piece on the curve, ${\rm f} (x)$, per definition has a less steep slope than the green piece ${\rm g} (x) $, we know that they form a concave hull, as illustrated by Figure \ref{fig:concave_hull_h2_value}. Hence, the hydrogen value curve can also be defined as
\begin{align} \label{h_value_min}
    {\rm \gamma} (x) = \text{min}\{{\rm g} (x) , {\rm f} (x) \}, \quad \forall x.
\end{align}

We now define the ${\rm g} (x) $ and ${\rm f} (x) $ functions. Under the assumption that the electrolyzer will always utilize the available green signal, we know that ${\rm g} (x) $ intercepts in the origin. Further, its slope is found as the value of the hydrogen produced for $x$ amount of energy.
\begin{align}
    {\rm g} (x) = \eta (\pi^{\rm gray} + \pi^{\rm green})x
\end{align}
The gray piece, however, does not intercept in origin unless there is no green-signal available. However, we know its slope and its value in $P^{\rm RES}$, which allows us to solve for its intercept $B$.
\begin{subequations}
    \begin{align}
    &{\rm f} (x) = \eta \pi^{\rm gray} x + B \\
    &{\rm f} (P^{\rm RES}) = {\rm g} (P^{\rm RES}) = \eta (\pi^{\rm gray} + \pi^{\rm green}) P^{\rm RES} \\
    &{\rm f} (P^{\rm RES}) = \eta (\pi^{\rm gray} + \pi^{\rm green}) P^{\rm RES} = \eta \pi^{\rm gray} P^{\rm RES} + B \\
    & B = \eta \pi^{\rm green} P^{\rm RES}
\end{align}
\end{subequations}

We achieve an expression for ${\rm f} (x)$ which depends on $P^{\rm RES}_s$ for its intercept.
\begin{align}
    &{\rm f} (x) = \eta \pi^{\rm gray} x + \eta \pi^{\rm green} P^{\rm RES}
\end{align}

The derived function for the hydrogen value function depends on $P^{\rm RES}$, which is uncertain. We draw a scenario of $P^{\rm RES}_s$, which is different in each scenario $s$. We therefore define the scenario dependent function as follows:
\begin{align}
    {\rm \gamma}_s (x) = \begin{cases}
        {\rm g}_s (x) = \eta (\pi^{\rm gray} + \pi^{\rm green})x  & \text{if} \quad x \leq P^{\rm RES}_s \\
        {\rm f}_s (x) = \eta \pi^{\rm gray} x + \eta \pi^{\rm green} P^{\rm RES}_s  & \text{if} \quad P^{\rm RES}_s < x \leq P^{\rm h}
    \end{cases}
\end{align}

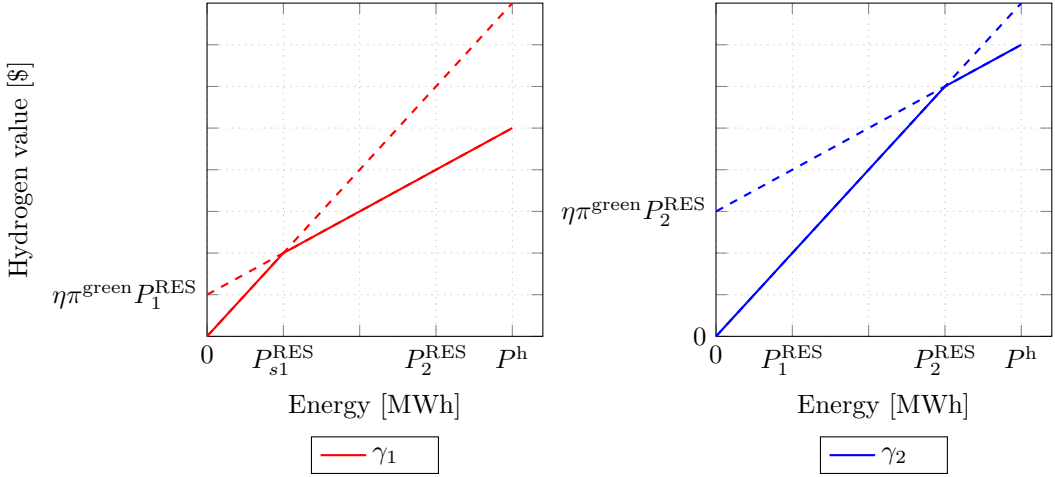
\begin{figure}[!ht]
\centering
\begin{tikzpicture}
    \begin{axis}[
        width=6cm,
        height=6cm,
        xlabel={Energy [MWh]},
        ylabel={Hydrogen value [\$]},
        grid=both,
        grid style={dotted},
        legend style={at={(0.5,-0.3)},anchor=north,legend columns=-1},
        ymin=0, ymax=200,
        xmin=0, xmax=110,
        xtick={0, 25, 50, 75, 100},
        ytick={0, 25, 50, 75, 100, 125, 150, 175 },
        yticklabels = {,$\eta\pi^{\rm green} P^{\rm RES}_{1}$,,,,,,},
        xticklabels = {0, $P^{\rm RES}_{s1}$, , $P^{\rm RES}_{2}$, $P^{\rm h}$}
    ]

        \addplot[forget plot, color=red, thick] coordinates {(0,0) (25,25*2)};
        \addplot[forget plot, dashed, color=red, thick] coordinates {(0,0) (100,100*2)};
        \addplot[color=red, thick] coordinates {(25,25*2) (100,75*1+25*2)};
        \addlegendentry{$\gamma_{1}\quad$};
        \addplot[forget plot, dashed, color=red, thick] coordinates {(0,25) (100,75*1+25*2)};
        \addlegendentry{$\gamma_{1}\quad$};
    \end{axis}
\end{tikzpicture}
\begin{tikzpicture}
    \begin{axis}[
        width=6cm,
        height=6cm,
        xlabel={Energy [MWh]},
        ylabel={},
        grid=both,
        grid style={dotted},
        legend style={at={(0.5,-0.3)},anchor=north,legend columns=-1},
        ymin=0, ymax=200,
        xmin=0, xmax=110,
        xtick={0, 25, 50, 75, 100},
        ytick={0, 25, 50, 75, 100, 125, 150, 175 },
        yticklabels = {0, , , $\eta\pi^{\rm green} P^{\rm RES}_{2}$,},
        xticklabels = {0, $P^{\rm RES}_{1}$, , $P^{\rm RES}_{2}$, $P^{\rm h}$}
    ]

        \addplot[forget plot, color=blue, thick] coordinates {(0,0) (75,75*2)};
        \addplot[forget plot, dashed, color=blue, thick] coordinates {(0,0) (100,100*2)};
        \addplot[color=blue, thick] coordinates {(75,75*2) (100,25*1+75*2)};
        \addlegendentry{$\gamma_{2}\quad$};
        \addplot[forget plot, dashed, color=blue, thick] coordinates {(0,75) (100,25*1+75*2)};
        \addlegendentry{$\gamma_{2}\quad$};
        
    \end{axis}
\end{tikzpicture}
\caption{Value of hydrogen for a given $P^{\rm RES}_s$.}
\label{hydrogen_value_two_s}
\end{figure}

Figure \ref{hydrogen_value_two_s} illustrates the hydrogen value as functions of the day-ahead purchase quantity for two different scenarios of the green-signal. The expected value of hydrogen production can be found as the sum of the probability weighed hydrogen value of all scenarios:
    \begin{align}
    \mathbb{E}[{\rm \gamma} (x)] = \sum_{s \in S} \rho_s {\rm \gamma}_s (x),
\end{align}
where $\rho_s$ is the probability of a scenario. 

\subsection{Formulating a profit-maximization problem for the electrolyzer under time-matching requirements}

In order to build its optimal bid curve, the electrolyzer should compute, for each possible values of the DA-price $\lambda^{\rm DA}$, the optimal quantity of electricity to buy $\boldsymbol{P^{\rm DA}}$ that maximizes its expected profit over all scenarios $s \in \mathcal{S}$ of the green-signal. Although the expected value of hydrogen expressed above is non-linear, we introduce an equivalent linear formulation of the profit-maximization problem of the electrolyzer for a given DA-price $\lambda^{\rm DA}$.

\subsubsection{Linear program formulation}

For a given value of the DA-price, we introduce the following linear program:
\begin{subequations}
\label{opt_problem}
\begin{align}
    & \max_{\boldsymbol{ p^{\rm DA}} , \boldsymbol{\gamma_{s\in \mathcal{S}}}} &&   Z=   \sum_{s \in \mathcal{S}} \big[ \rho_s \boldsymbol{\gamma_{s} } \big]  - \lambda^{\rm DA} \boldsymbol{ p^{\rm DA}}  \\
    & && \boldsymbol{\gamma_{s} }  \leq  \eta (\pi^{\rm gray} + \pi^{\rm green})  \boldsymbol{ p^{\rm DA}} &&   :\mu^{\rm green}_s && \forall s \in \mathcal{S} \label{ineq1} \\
    & && \boldsymbol{\gamma_{s} }  \leq  \eta \pi^{\rm gray}  \boldsymbol{ p^{\rm DA}} + \eta  \pi^{\rm green}  P^{\rm RES}_s &&   :\mu^{\rm gray}_s && \forall s \in \mathcal{S}  \label{ineq2} \\
    & && - \boldsymbol{p^{\rm DA} } \leq 0 && : \mu^{\rm lb} \\
    & && \boldsymbol{p^{\rm DA} } \leq P^{\rm h} && : \mu^{\rm up},
\end{align}
\end{subequations}
where $P^{\rm h}$ represents the load of the electrolyzer, $\eta$ its constant efficiency, $\pi^{\rm gray}$ the gray-hydrogen price, and $\pi^{\rm green}$ the green-hydrogen bonus. In order to discard trivial cases, we assume that all these parameters are strictly positive.

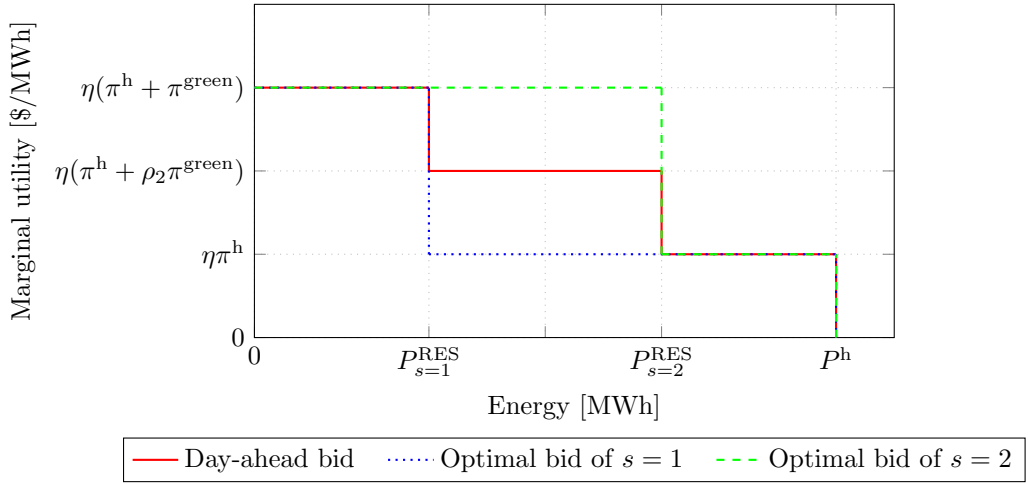
\begin{figure}[!ht]
\centering
\begin{tikzpicture}
    \begin{axis}[
        width=10cm,
        height=6cm,
        xlabel={Energy [MWh]},
        ylabel={Marginal utility [\$/MWh]},
        grid=both,
        grid style={dotted},
        legend style={at={(0.5,-0.3)},anchor=north,legend columns=-1},
        ymin=0, ymax=80,
        xmin=0, xmax=110,
        xtick={0, 30, 50, 70, 100},
        ytick={0, 20, 20+40/2, 60},
        yticklabels = {0, $\eta\pi^{\rm h}$, $\eta(\pi^{\rm h}+ \rho_2 \pi^{\rm green})$, $\eta(\pi^{\rm h}+\pi^{\rm green})$,  $\eta(\pi^{\rm h}+\pi^{\rm green})$ },
        xticklabels = {0, $P^{\rm RES}_{s = 1}$, , $P^{\rm RES}_{s = 2}$, $P^{\rm h}$} 
    ]


        \addplot[forget plot, color=red, thick, ] coordinates {(0,60) (30,60)};
        
        \addplot[forget plot, color=red, thick] coordinates {(30,60) (30,40)};

        \addplot[forget plot, color=red, thick] coordinates {(30,40) (70,40)};

        \addplot[forget plot, color=red, thick] coordinates {(70,40) (70,20)};

        \addplot[forget plot, color=red, thick] coordinates {(100,20) (70,20)};

        \addplot[color=red, thick] coordinates {(100,0) (100,20)};
        \addlegendentry{Day-ahead bid$\quad$};
        
        \addplot[forget plot, dotted, color=blue, thick] coordinates {(0,60) (30,60)};
        
        \addplot[forget plot, dotted, color=blue, thick] coordinates {(30,20) (30,60)};

        \addplot[forget plot, dotted, color=blue, thick] coordinates {(30,20) (100,20)};

        \addplot[dotted, color=blue, thick] coordinates {(100,0) (100,20)};
        \addlegendentry{Optimal bid of $s = 1\quad$};

        \addplot[forget plot, dashed, color=green, thick] coordinates {(0,60) (70,60)};
        
        \addplot[forget plot, dashed, color=green, thick] coordinates {(70,20) (70,60)};

        \addplot[forget plot, dashed, color=green, thick] coordinates {(70,20) (100,20)};

        \addplot[dashed, color=green, thick] coordinates {(100,0) (100,20)};
        \addlegendentry{Optimal bid of $s = 2$};

    \end{axis}
\end{tikzpicture}
\caption{Optimal electrolyzer bid curve under two (increasing) scenarios of $P^{\rm RES}$.}
\label{bid_curve_illustration}
\end{figure}

Figure \ref{bid_curve_illustration} illustrates the optimal quantity $\boldsymbol{p^{\rm DA}}$ (in red), obtained from solving the linear program \eqref{opt_problem}, over a range of DA-prices and for two scenarios of the green-signal.

\subsubsection{KKT conditions}

We derive the Lagrangian function of \eqref{opt_problem}:
\begin{subequations}
    \begin{align}
    \mathcal L(p^{\rm DA}, \gamma_{s \in \mathcal{S}}, \mu^{\rm (\cdot)}) = &  \lambda^{\rm DA} \boldsymbol{ p^{\rm DA}} - \sum_{s \in \mathcal{S}} \big[ \rho_s \boldsymbol{\gamma_{s} } \big]    - \mu^{\rm lb} \boldsymbol{p^{\rm DA} } + \mu^{\rm ub} (\boldsymbol{p^{\rm DA} } - P^{\rm h} ) \notag \\ 
    &+ \sum_{s \in \mathcal{S}} \big[ \mu^{\rm green}_s \left(\boldsymbol{\gamma_{s} }  - \eta (\pi^{\rm gray} + \pi^{\rm green})  \boldsymbol{ p^{\rm DA}} \right)  \\ 
    & +  \mu^{\rm gray}_s \left( \boldsymbol{\gamma_{s} }  -  \eta \pi^{\rm gray}  \boldsymbol{ p^{\rm DA}} - \eta  \pi^{\rm green}  P^{\rm RES}_s \right) \big], \notag
\end{align}
\end{subequations}
and the following KKT conditions: \\
\noindent \textit{1) Stationarity conditions:}
\begin{subequations}
    \begin{align}
&\frac{\partial \mathcal L(.)}{\partial p^{\rm  DA}} = \lambda^{\rm DA} -\mu^{\rm lb} + \mu^{\rm ub} + \sum_{s \in \mathcal{S}} \big[ - \mu^{\rm green}_s \eta (\pi^{\rm gray} + \pi^{\rm green}) - \mu^{\rm gray}_s \eta \pi^{\rm gray} \big] = 0 , \label{stat_cond_1} \\
& \frac{\partial \mathcal L(.)}{\partial \gamma_{s}} = - \rho_s + \mu^{\rm green}_s +  \mu^{\rm gray}_s =  0 \qquad \forall s \in \mathcal{S}. \label{stat_cond_2}
\end{align}
\noindent \textit{2) Complementary slackness conditions:}
\begin{align} \label{comp_cond}
& (\boldsymbol{\gamma_{s} } - \eta (\pi^{\rm gray} + \pi^{\rm green})  \boldsymbol{ p^{\rm DA}} ) \mu^{\rm green}_s = 0  &&\forall s \in \mathcal{S}, \\
& (\boldsymbol{\gamma_{s} }  -  \eta \pi^{\rm gray}  \boldsymbol{ p^{\rm DA}} - \eta  \pi^{\rm green}  P^{\rm RES}_s) \mu^{\rm gray}_s = 0  && \forall s \in \mathcal{S}, \\
& - \boldsymbol{p^{\rm DA} }  \mu^{\rm lb}  = 0 , \\
 & (\boldsymbol{p^{\rm DA} } - P^{\rm h}) \mu^{\rm ub} = 0.
\end{align}
\noindent \textit{3) Primal and dual feasibility conditions:}
\begin{align}
     &  \boldsymbol{\gamma_{s} }  \leq  \eta (\pi^{\rm gray} + \pi^{\rm green})  \boldsymbol{ p^{\rm DA}} && \forall s \in \mathcal{S}, \\
    &  \boldsymbol{\gamma_{s} }  \leq  \eta \pi^{\rm gray}  \boldsymbol{ p^{\rm DA}} + \eta  \pi^{\rm green}  P^{\rm RES}_s && \forall s \in \mathcal{S}, \\
    &  \boldsymbol{p^{\rm DA} } \leq P^{\rm h}, \\
    & \boldsymbol{ p^{\rm DA}} , \boldsymbol{\gamma_{s\in \mathcal{S}}}, \mu^{(\cdot)} \geq 0. \label{KKT_last}
\end{align}
\end{subequations}

\subsubsection{Equivalence to profit-maximization problem of electrolyzer}

We can show that the linear program \eqref{opt_problem} is equivalent to the profit-maximization problem of the electrolyzer for a given DA-price. Firstly, we observe that the first-order stationarity conditions \eqref{stat_cond_2}, ensure that, at optimality, for each scenario $s \in \mathcal{S}$, at least one of the dual variables $\mu_s^{green}$, $\mu_s^{green}$ must be strictly positive, and by extension at least one of the inequality constraints \eqref{ineq1}-\eqref{ineq2} must be binding. This ensures that 

$$\boldsymbol{\gamma_s} = \min \left(\eta (\pi^{\rm gray} + \pi^{\rm green})  \boldsymbol{ p^{\rm DA}},\eta \pi^{\rm gray}  \boldsymbol{ p^{\rm DA}} + \eta  \pi^{\rm green}  P^{\rm RES}_s \right) = \rm \gamma_s (\boldsymbol{ p^{\rm DA}}), \quad \forall s \in \mathcal{S}. $$
It follows from the definition in \eqref{h_value_min} that the expected value of hydrogen $\mathbb{E}[{\rm \gamma} (\boldsymbol{ p^{\rm DA}})] = \sum_{s \in \mathcal{S}} \big[ \rho_s \boldsymbol{\gamma_s} \big]$ for a given amount of day-ahead purchase quantity $\boldsymbol{ p^{\rm DA}}$. As a result, the linear program \eqref{opt_problem} is equivalent to the optimization problem of an electrolyzer aiming at finding the optimal bid quantity $\boldsymbol{p^{\rm DA}}$ that maximizes its expected profit over all scenarios of the green signal, for a given DA-price.

\subsection{Interpretation of KKT conditions} 

For a given DA-price $\lambda^{\rm DA}$, the complementary slackness conditions \eqref{comp_cond} can be interpreted as the following combinations of conditions that all optimal solutions to program \eqref{opt_problem} must satisfy:
\begin{subequations}
\begin{align}
    & (\boldsymbol{\gamma_{s} } - \eta (\pi^{\rm gray} + \pi^{\rm green})  \boldsymbol{ p^{\rm DA}} ) = 0 & \text{and-or} && \mu^{\rm green}_s = 0 && \forall s \in \mathcal{S} \label{comp_s1}\\
    & (\boldsymbol{\gamma_{s} }  -  \eta \pi^{\rm gray}  \boldsymbol{ p^{\rm DA}} - \eta  \pi^{\rm green}  P^{\rm RES}_s)  = 0 &\text{and-or} &&\mu^{\rm gray}_s = 0 && \forall s \in \mathcal{S} \label{comp_s2} \\
    & \boldsymbol{ p^{\rm DA}} = 0 & \text{and-or} && \mu^{\rm lb} = 0 \label{comp_DA_lb}\\
    &(\boldsymbol{ p^{\rm DA}} - P^{\rm h}) = 0 &\text{and-or} &&\mu^{\rm ub} = 0 \label{comp_DA_ub}
\end{align}
\end{subequations}

As \eqref{comp_DA_lb}-\eqref{comp_DA_ub} are scenario independent and solely depend on the optimal quantity decision $p^{\rm DA}$, they can be intuitively interpreted and provide 3 distinct cases of optimal day-ahead purchase quantity:
\begin{enumerate}
    \item \textit{No-bid case:} If the optimal solution is $\boldsymbol{p^{\rm DA}} = 0$, then $\boldsymbol{p^{\rm DA}} < P^{\rm h}$.\footnote{It is assumed that $P^{\rm h}>0$.} It follows that $\mu^{\rm ub} = 0$.
    \item \textit{Partial-capacity bid case:} If $0 < p^{\rm DA} < P^{\rm h}$, then, both $\mu^{\rm ub} = \mu^{\rm lb} = 0$.
    \item \textit{Maximum-capacity bid case:} By analogy, if $p^{\rm DA} = P^{\rm h}$,  then $\mu^{\rm lb} = 0$.
\end{enumerate}  

In addition, the complementary slackness conditions \eqref{comp_s1}-\eqref{comp_s2} and the stationarity conditions \eqref{stat_cond_2} provide combinations of conditions on the scenarios of hydrogen values $\boldsymbol{\gamma_s}$, that depend on the relative value of the green signal scenarios $P_s^{\rm RES}$ and the optimal day-ahead purchase quantity $\boldsymbol{P^{\rm DA}}$. Therefore, for a given solution $\boldsymbol{P^{\rm DA}}$, we introduce the following subsets of scenarios:

\begin{align}
    &\mathcal{S}^{\rm green}(\boldsymbol{P^{\rm DA}}) := \{s \in \mathcal{S} \mid \boldsymbol{P^{\rm DA}} < P^{\rm RES}_s \}, \\
    &\mathcal{S}^{\rm RES}(\boldsymbol{P^{\rm DA}}) := \{s \in \mathcal{S} \mid \boldsymbol{P^{\rm DA}} = P^{\rm RES}_s \}, \\
    &\mathcal{S}^{\rm gray}(\boldsymbol{P^{\rm DA}}) := \{s \in \mathcal{S} \mid \boldsymbol{P^{\rm DA}} > P^{\rm RES}_s \}. 
\end{align}

By definition, these subsets form a partition of $
\mathcal{S}$, i.e. $\mathcal{S}^{\rm green}(\boldsymbol{P^{\rm DA}})\cup\mathcal{S}^{\rm RES}(\boldsymbol{P^{\rm DA}})\cup\mathcal{S}^{\rm gray}(\boldsymbol{P^{\rm DA}})=\mathcal{S}$, and are pairwise disjoint, i.e., $\mathcal{S}^{\rm green}(\boldsymbol{P^{\rm DA}})\cap\mathcal{S}^{\rm RES}(\boldsymbol{P^{\rm DA}})=\emptyset$, $\mathcal{S}^{\rm green}(\boldsymbol{P^{\rm DA}})\cap\mathcal{S}^{\rm gray}(\boldsymbol{P^{\rm DA}}) =  \emptyset$, and $\mathcal{S}^{\rm RES}(\boldsymbol{P^{\rm DA}})\cap\mathcal{S}^{\rm gray}(\boldsymbol{P^{\rm DA}}) =  \emptyset$.

We can now derive the following optimality conditions, for a given value of the optimal day-ahead purchase quantity $\boldsymbol{P^{\rm DA}}$:
\begin{enumerate}[label=\roman*]
    \item For all $s \in \mathcal{S}^{\rm green}(\boldsymbol{P^{\rm DA}})$, $\boldsymbol{\gamma_{s} } = \eta (\pi^{\rm gray} + \pi^{\rm green}) \boldsymbol{p^{\rm DA}} < \eta (\pi^{\rm gray} \boldsymbol{p^{\rm DA}} + \pi^{\rm green}) P_s^{\rm RES}$. It follows that $\mu^{\rm gray}_s = 0$ and $\mu^{\rm green}_s = \rho_s >0$.
    \item For all $s \in \mathcal{S}^{\rm RES}(\boldsymbol{P^{\rm DA}})$, then $\boldsymbol{\gamma_{s} } = \eta (\pi^{\rm gray} + \pi^{\rm green}) \boldsymbol{p^{\rm DA}} = \eta (\pi^{\rm gray} \boldsymbol{p^{\rm DA}} + \pi^{\rm green}) P_s^{\rm RES}$. It follows that both $\mu^{\rm gray}_s \geq 0$ and $\mu^{\rm gray}_s \geq 0$, and $\mu^{\rm gray}_s + \mu^{\rm gray}_s = \rho_s$.
    \item For all $s \in \mathcal{S}^{\rm gray}(\boldsymbol{P^{\rm DA}})$, then $\boldsymbol{\gamma_{s} } = \eta (\pi^{\rm gray} \boldsymbol{p^{\rm DA}} + \pi^{\rm green}) P_s^{\rm RES} < \eta (\pi^{\rm gray} + \pi^{\rm green}) \boldsymbol{p^{\rm DA}}$. It follows that $\mu^{\rm green}_s = 0$ and $\mu^{\rm gray}_s = \rho_s >0$.
\end{enumerate}

\color{black} 

\subsection{Conditions linking DA-prices and optimal day-ahead purchase quantity}

We can now derive two types of conditions linking the day-ahead prices $\lambda^{\rm DA}$ and optimal day-ahead purchase quantities $\boldsymbol{P^{\rm DA}}$, namely: 
\begin{itemize}
    \item The sets $\mathcal{Q}(\lambda^{\rm DA}) := \{\boldsymbol{P^{\rm DA}} \in \mathbb{R} \mid  \boldsymbol{P^{\rm DA}} \text{ optimal solution of \eqref{opt_problem} under the value } \lambda^{\rm DA} \}$,  representing the range of optimal day-ahead purchase quantities $\boldsymbol{P^{\rm DA}}$ for a given day-ahead price $\lambda^{\rm DA}$. The notation of $\mathcal{Q}$ indicates the set of optimal \textit{quantities}.
    \item The sets $\mathcal{P}(\boldsymbol{P^{\rm DA}}) := \{ \lambda^{\rm DA} \in \mathbb{R} \mid  \boldsymbol{P^{\rm DA}} \text{ optimal solution of \eqref{opt_problem} under the value } \lambda^{\rm DA} \}$, representing the range of day-ahead prices $\lambda^{\rm DA}$ for which a given day-ahead purchase quantity $\boldsymbol{P^{\rm DA}}$ is optimal. The notation of $\mathcal{P}$ indicates the set of optimal \textit{prices}.
\end{itemize}
Generally, these conditions mean that a range of prices can be optimal for a given quantity, or that a range of quantities are optimal for a given price. These conditions will provide the basis for the derivation of the electrolyzer's optimal price-quantity bid curve, each corresponding a specific segment of the curve as further discussed in Section \ref{appendix_bid_curve}. We derive below these conditions for each case of optimal day-ahead purchase quantity, i.e. the \textit{no-bid}, \textit{partial-capacity bid},  and \textit{maximum-capacity bid} cases.

\subsubsection{No-bid ($\mathbf{p^{\rm \textbf{DA}} = 0}$)}

We recall that, at optimality, $\mu^{\rm ub} = 0$. As a result, the stationarity condition \eqref{stat_cond_1} becomes:
\begin{subequations}
\begin{align} \label{stad_cond_no}
    & \lambda^{\rm DA} -\mu^{\rm lb} + \sum_{s \in \mathcal{S}} \big[ - \mu^{\rm green}_s \eta (\pi^{\rm gray} + \pi^{\rm green}) - \mu^{\rm gray}_s \eta \pi^{\rm gray} \big] = 0
\end{align}
Following the definition of the subsets of scenarios $\mathcal{S}^{\rm green}(\boldsymbol{P^{\rm DA}})$, $\mathcal{S}^{\rm RES}(\boldsymbol{P^{\rm DA}})$, and $\mathcal{S}^{\rm gray}(\boldsymbol{P^{\rm DA}})$, the values of $\mu_s^{\rm green}$ and $\mu_s^{\rm gray}$ can be substituted in \eqref{stad_cond_no}, leading to: 
\begin{align}
    &\lambda^{\rm DA} - \mu^{\rm lb}- \sum_{\mathclap{s \in \mathcal{S}^{\rm green}(\boldsymbol{P^{\rm DA}})}} \big[ \rho_s \eta (\pi^{\rm gray} + \pi^{\rm green})  \big] - \sum_{\mathclap{s \in \mathcal{S}^{\rm gray}(\boldsymbol{P^{\rm DA}})}} \big[ \rho_s \eta \pi^{\rm gray} \big] -\sum_{\mathclap{s \in \mathcal{S}^{\rm RES}(\boldsymbol{P^{\rm DA}})}} \big[ \rho_s  \eta \pi^{\rm gray}  + \mu^{\rm green}_s \eta \pi^{\rm green} \big]  = 0 \\
    &\lambda^{\rm DA} - \mu^{\rm lb} - \eta \pi^{\rm gray} \sum_{s \in \mathcal{S}}  \rho_s   - \eta  \pi^{\rm green} \sum_{\mathclap{s \in \mathcal{S}^{\rm green}(\boldsymbol{P^{\rm DA}})}}  \rho_s   - \eta \pi^{\rm green} \sum_{\mathclap{s \in \mathcal{S}^{\rm RES}(\boldsymbol{P^{\rm DA}})}}   \mu^{\rm green}_s   = 0 \\
    &\lambda^{\rm DA} - \eta \pi^{\rm gray}   - \eta  \pi^{\rm green} \sum_{\mathclap{s \in \mathcal{S}^{\rm green}(\boldsymbol{P^{\rm DA}})}}  \rho_s   - \eta \pi^{\rm green} \sum_{\mathclap{s \in \mathcal{S}^{\rm RES}(\boldsymbol{P^{\rm DA}})}}   \mu^{\rm green}_s   = \mu^{\rm lb} \geq 0 \\
    &\lambda^{\rm DA} - \eta \pi^{\rm gray}   - \eta  \pi^{\rm green} \sum_{\mathclap{s \in \mathcal{S}^{\rm green}(\boldsymbol{P^{\rm DA}})}}  \rho_s   \geq  \eta \pi^{\rm green} \sum_{\mathclap{s \in \mathcal{S}^{\rm RES}(\boldsymbol{P^{\rm DA}})}}   \mu^{\rm green}_s 
\end{align}
We derive the following bounds on the right-hand-side of this equation from the stationarity condition \eqref{stat_cond_2}:
\begin{align}
        0 \leq \eta \pi^{\rm green} \sum_{\mathclap{s \in \mathcal{S}^{\rm RES}(\boldsymbol{P^{\rm DA}})}}  \mu^{\rm green}_s  \leq \eta \pi^{\rm green} \sum_{\mathclap{s \in \mathcal{S}^{\rm RES}(\boldsymbol{P^{\rm DA}})}} \rho_s,
\end{align}
leading to the following lower-bound on the day-ahead prices:
\begin{align}
    \lambda^{\rm DA} \geq  \eta \pi^{\rm gray}   + \eta  \pi^{\rm green} \sum_{\mathclap{s \in \mathcal{S}^{\rm green}(\boldsymbol{P^{\rm DA}})}}  \rho_s.
\end{align}

Given that $P^{\rm DA} = 0 \leq P^{\rm RES}_s$, $\forall s \in \mathcal{S}$, we can directly formulate the following range of day-ahead prices for which $\boldsymbol{P^{\rm DA}} = 0$ is an optimal solution as:
\begin{align}
    \mathcal{P}(\boldsymbol{P^{\rm DA}}=0) = \Big[ \eta \pi^{\rm gray}   + \eta  \pi^{\rm green} \sum_{\mathclap{s \in \mathcal{S}^{\rm green}(\boldsymbol{P^{\rm DA}}=0)}} \rho_s \quad , \ +\infty \Big[ ,
\end{align}
where $\mathcal{S}^{\rm green}(\boldsymbol{P^{\rm DA}}=0) = \mathcal{S}\setminus\{1\}$ if $P_1^{\rm RES} = 0$ and $\mathcal{S}^{\rm green}(\boldsymbol{P^{\rm DA}}=0) = \mathcal{S}$ if $P^{\rm RES} > 0$\footnote{Indeed, due to the definition of the unique and increasing scenarios, a scenario $s$ can only belong to $\mathcal{S}^{\rm RES}(\boldsymbol{P^{\rm DA}} = 0)$ if $s=1$ and $P_1^{\rm RES} = 0$.}. In other words, these ranges can be defined explicitly as:
\begin{align}
    & \mathcal{P}(\boldsymbol{P^{\rm DA}}=0) = \begin{cases} & \big[ \eta \pi^{\rm gray}   + \eta  \pi^{\rm green} (1 - \rho_1) , \ +\infty \big[ , \text{ if }P^{\rm RES}_1  = 0 \\
    &  \big[ \eta \pi^{\rm gray}   + \eta  \pi^{\rm green}  , \ +\infty \big [ , \text{ if }P^{\rm RES}_1 >0.
    \end{cases}
\end{align}

\end{subequations}

\subsubsection{Partial-capacity bid ($\mathbf{0 < p^{\rm \textbf{DA}} < P^{\rm \textbf{h}}}$)}

We recall that, at optimality, $\mu^{\rm lb} = \mu^{\rm ub} = 0$. As a result the stationarity condition \eqref{stat_cond_1} becomes:
\begin{subequations}
\begin{align} \label{stad_cond_partial}
    & \lambda^{\rm DA} + \sum_{\mathclap{s \in \mathcal{S}}} \big[ - \mu^{\rm green}_s \eta (\pi^{\rm gray} + \pi^{\rm green}) - \mu^{\rm gray}_s \eta \pi^{\rm gray} \big] = 0
\end{align}
Following the definition of the subsets of scenarios $\mathcal{S}^{\rm green}(\boldsymbol{P^{\rm DA}})$, $\mathcal{S}^{\rm RES}(\boldsymbol{P^{\rm DA}})$, and $\mathcal{S}^{\rm gray}(\boldsymbol{P^{\rm DA}})$, the values of $\mu_s^{\rm green}$ and $\mu_s^{\rm gray}$ can be substituted in \eqref{stad_cond_partial}, leading to: 
\begin{align}
    & \lambda^{\rm DA} -  \sum_{\mathclap{s \in \mathcal{S}^{\rm green}(\boldsymbol{P^{\rm DA}})}} \big[ \rho_s \eta (\pi^{\rm gray} + \pi^{\rm green})  \big] - \sum_{\mathclap{s \in \mathcal{S}^{\rm gray}(\boldsymbol{P^{\rm DA}})}} \big[ \rho_s \eta \pi^{\rm gray} \big] -\sum_{\mathclap{s \in \mathcal{S}^{\rm RES}(\boldsymbol{P^{\rm DA}})}} \big[ \rho_s  \eta \pi^{\rm gray}  + \mu^{\rm green}_s \eta \pi^{\rm green} \big]  = 0 \\
    & \lambda^{\rm DA} - \sum_{{s \in \mathcal{S}}} \big[ \rho_s \eta \pi^{\rm gray} \big] -\sum_{\mathclap{s \in \mathcal{S}^{\rm green}(\boldsymbol{P^{\rm DA}})}} \big[ \rho_s \eta  \pi^{\rm green}  \big]  -\sum_{\mathclap{s \in \mathcal{S}^{\rm RES}(\boldsymbol{P^{\rm DA}})}} \big[   \mu^{\rm green}_s \eta \pi^{\rm green} \big]  = 0 \\
    & \lambda^{\rm DA} - \eta \pi^{\rm gray}  -\eta \pi^{\rm green} \sum_{\mathclap{s \in \mathcal{S}^{\rm green}(\boldsymbol{P^{\rm DA}})}} \rho_s    -\eta \pi^{\rm green} \sum_{\mathclap{s \in \mathcal{S}^{\rm RES}(\boldsymbol{P^{\rm DA}})}} \mu^{\rm green}_s  = 0 \\
    & \lambda^{\rm DA} - \eta \pi^{\rm gray}  -\eta  \pi^{\rm green} \sum_{\mathclap{s \in \mathcal{S}^{\rm green}(\boldsymbol{P^{\rm DA}})}} \rho_s = \eta \pi^{\rm green} \sum_{\mathclap{s \in \mathcal{S}^{\rm RES}(\boldsymbol{P^{\rm DA}})}} \mu^{\rm green}_s
\end{align}
We derive the following bounds on the right-hand-side of this equation, from the stationarity condition \eqref{stat_cond_2}:
\begin{align}
        0 \leq \eta \pi^{\rm green} \sum_{\mathclap{s \in \mathcal{S}^{\rm RES}(\boldsymbol{P^{\rm DA}})}} \mu^{\rm green}_s \leq \eta \pi^{\rm green} \sum_{\mathclap{s \in \mathcal{S}^{\rm RES}(\boldsymbol{P^{\rm DA}})}}\rho_s ,
\end{align}
leading to the following upper- and lower-bounds on the range of DA-prices for which the day-ahead purchase quantity $\boldsymbol{P^{\rm DA}}$ is an optimal solution:
\begin{align} \label{bounds_DA_price_PCB}
        & \eta \pi^{\rm gray}  + \eta  \pi^{\rm green} \sum_{\mathclap{s \in \mathcal{S}^{\rm green}(\boldsymbol{P^{\rm DA}})}} \rho_s \leq \lambda^{\rm DA} \leq \eta \pi^{\rm gray} + \eta  \pi^{\rm green} \sum_{\mathclap{s \in \mathcal{S}^{\rm green}(\boldsymbol{P^{\rm DA}})}} \rho_s + \eta \pi^{\rm green} \sum_{\mathclap{s \in \mathcal{S}^{\rm RES}(\boldsymbol{P^{\rm DA}})}}\rho_s.
\end{align}

To derive the conditions linking the optimal price and quantity bids for the electrolyzer, we must consider all possible configurations of the day-ahead purchase quantity $0 < \boldsymbol{P^{\rm DA}} < P^{\rm}$ relative to the green signal scenarios $P^{\rm RES}_s$, for $s \in \mathcal{S}$. Therefore, we study the following cases, based on the values of the subsets $\mathcal{S}^{\rm RES}(\boldsymbol{P^{\rm DA}})$ and $\mathcal{S}^{\rm green}(\boldsymbol{P^{\rm DA}})$:
\begin{enumerate}[label=\alph*.]
    \item If $\boldsymbol{P^{\rm DA}} \notin \{P^{\rm RES}_s : s \in \mathcal{S} \}$ (i.e. $\mathcal{S}^{\rm RES}(\boldsymbol{P^{\rm DA}}) = \emptyset$): In this case, following the definition of the scenarios of the green signal, there are (up-to) $s+1$ possible ranges of values for the day-ahead purchase quantity $\boldsymbol{P^{\rm DA}}$, namely:
    $\mathcal{Q}_1=]0,P^{\rm RES}_1[$ , $\mathcal{Q}_s = ]P^{\rm RES}_{s-1},P^{\rm RES}_{s}[$ for $s \in \mathcal{S}\setminus\{1\}$, and $\mathcal{Q}_{S+1} = ]P^{\rm RES}_S,P^{\rm h}[$. Indeed, by definition, the ranges $]0,P^{\rm RES}_1[$ and $]P^{\rm RES}_S,P^{\rm h}[$ are empty if $P^{\rm RES}_1=0$ or $P^{\rm RES}_S = P^{\rm h}$. 
It follows from \eqref{bounds_DA_price_PCB}, that for each value of the day-ahead price $$ \lambda^{\rm DA} = \eta \pi^{\rm gray}  + \eta  \pi^{\rm green} \sum_{\mathclap{i \in \mathcal{S}^{\rm green}(\boldsymbol{P^{\rm DA}} \in \mathcal{Q}_s)}} \rho_i \quad , $$ all day-ahead purchase quantities $\boldsymbol{P^{\rm DA}} \in \mathcal{Q}_s$ are optimal. That is,
    \begin{align}
    \mathcal{Q} \big( \lambda^{\rm DA} = \eta \pi^{\rm gray}  + \eta  \pi^{\rm green} \sum_{\mathclap{i \in \mathcal{S}^{\rm green}(\boldsymbol{P^{\rm DA}} \in \mathcal{Q}_s)}} \rho_i \quad \ \ \ \big) = \mathcal{Q}_s, \quad \forall s \in \mathcal{S}.
    \end{align}
    In addition, we observe that, over each range $\mathcal{Q}_s$, the corresponding subset of green scenarios is given by $\mathcal{S}^{\rm green}(\boldsymbol{P^{\rm DA}} \in \mathcal{Q}_s) = \{s,...,S\}$ and $\mathcal{S}^{\rm green}(\boldsymbol{P^{\rm DA}} \in \mathcal{Q}_{S+1}) = \emptyset$. Therefore, we can explicitly formulate these ranges as:
    \begin{align}
    & \mathcal{Q} ( \lambda^{\rm DA} = \eta \pi^{\rm gray}  + \eta  \pi^{\rm green}  ) = \big]0,P^{\rm RES}_1\big[,  \\
    & \mathcal{Q} ( \lambda^{\rm DA} = \eta \pi^{\rm gray}  + \eta  \pi^{\rm green} \sum_{\mathclap{i \in \{s,...S\} }} \rho_i ) = \big]P^{\rm RES}_{s-1},P^{\rm RES}_{s}\big[, \quad \forall s \in \{2,...,S\} \\
    & \mathcal{Q} ( \lambda^{\rm DA} = \eta \pi^{\rm gray} ) = \big]P^{\rm RES}_{S},P^{\rm h}\big[.
    \end{align}
    \item If $\boldsymbol{P^{\rm DA}} \in \{P^{\rm RES}_s : s \in \mathcal{S} \}$ (i.e. $\mathcal{S}^{\rm RES}(\boldsymbol{P^{\rm DA}}) \neq \emptyset$): In this case, the day-ahead purchase quantity $0 < \boldsymbol{P^{\rm DA}} < P^{\rm 0}$ can take (up-to) $s$ possible values $\boldsymbol{P^{\rm DA}}  = P_s^{\rm RES}$. Indeed, if $P_1^{\rm RES}=0$ and/or $P_s^{\rm RES}=P^{\rm h}$, these values are infeasible for the partial-capacity case where $0 < \boldsymbol{P^{\rm DA}} < P^{\rm 0}$. 
    And, for each possible $s \in \mathcal{S}$ s.t. $0 < P^{\rm RES} <P^{\rm h}$, the range of day-ahead prices $\lambda^{\rm DA}$ that make the day-ahead purchase quantity $\boldsymbol{P^{\rm DA}}  = P_s^{\rm RES}$ optimal is given by: 
    \begin{align}
        & \mathcal{P} (\boldsymbol{P^{\rm DA}}  = P_s^{\rm RES}) = \big[ \eta \pi^{\rm gray}  + \eta  \pi^{\rm green} \sum_{\mathclap{i \in \mathcal{S}^{\rm green}(\boldsymbol{P^{\rm DA}}=P^{\rm RES}_s)}} \rho_i \qquad \ , \eta \pi^{\rm gray} + \eta  \pi^{\rm green} \sum_{\mathclap{i \in \mathcal{S}^{\rm green}(\boldsymbol{P^{\rm DA}}=P^{\rm RES}_s)}} \rho_i \quad + \eta  \pi^{\rm green} \rho_s  \quad \big],
    \end{align}
    where the corresponding subsets of green scenarios are given by $\mathcal{S}^{\rm green}(\boldsymbol{P^{\rm DA}} = P^{\rm RES}_s) = \{s+1,...S\}$ for $s \in \mathcal{S}\setminus{S}$ and $\mathcal{S}^{\rm green}(\boldsymbol{P^{\rm DA}} = P^{\rm RES}_S) = \emptyset$. Therefore, we can explicitly formulate these ranges as: 
    \begin{align}
        & \mathcal{P} (\boldsymbol{P^{\rm DA}}  = P_1^{\rm RES}) = \big[ \eta \pi^{\rm gray}  + \eta  \pi^{\rm green} \sum_{\mathclap{i \in \{2,...,S\}}} \rho_i \ , \eta \pi^{\rm gray} + \eta  \pi^{\rm green} \sum_{\mathclap{i \in \{1,...,S\}}} \rho_i \big] , \text{ if } P_1^{\rm RES} > 0 \\
        & \mathcal{P} (\boldsymbol{P^{\rm DA}}  = P_s^{\rm RES}) = \big[ \eta \pi^{\rm gray}  + \eta  \pi^{\rm green} \sum_{\mathclap{i \in \{s+1,...,S\}}} \rho_i \ , \eta \pi^{\rm gray} + \eta  \pi^{\rm green} \sum_{\mathclap{i \in \{s,...,S\}}} \rho_i \big] , \forall s \in \{2,...,S-1\} \\
        & \mathcal{P} (\boldsymbol{P^{\rm DA}}  = P_{S}^{\rm RES}) = \big[ \eta \pi^{\rm gray}  \ , \eta \pi^{\rm gray} + \eta \pi^{\rm green} \rho_S \big] , \text{ if } P_S^{\rm RES} < P^{\rm h} .
    \end{align}

\color{black}
    \end{enumerate}

\end{subequations}

\subsubsection{Maximum-capacity bid ($\mathbf{p^{\rm \textbf{DA}} = P^{\rm \textbf{h}}}$)}

We recall that, at optimality $\mu^{\rm lb} = 0$. 
As a result the stationarity condition \eqref{stat_cond_1} becomes:
\begin{subequations}
\begin{align}
    & \lambda^{\rm DA} + \mu^{\rm ub} + \sum_{s \in \mathcal{S}} \big[ - \mu^{\rm green}_s \eta (\pi^{\rm gray} + \pi^{\rm green}) - \mu^{\rm gray}_s \eta \pi^{\rm gray} \big] = 0 
\end{align}
Following the definition of the subsets of scenarios $\mathcal{S}^{\rm green}(\boldsymbol{P^{\rm DA}})$, $\mathcal{S}^{\rm RES}(\boldsymbol{P^{\rm DA}})$, and $\mathcal{S}^{\rm gray}(\boldsymbol{P^{\rm DA}})$, the values of $\mu_s^{\rm green}$ and $\mu_s^{\rm gray}$ can be substituted in \eqref{stad_cond_no}, leading to:
\begin{align}
    &\lambda^{\rm DA} + \mu^{\rm ub} -  \sum_{\mathclap{s \in \mathcal{S}^{\rm green}(\boldsymbol{P^{\rm DA}})}} \big[ \rho_s \eta (\pi^{\rm gray} + \pi^{\rm green})  \big] - \sum_{\mathclap{s \in \mathcal{S}^{\rm gray}(\boldsymbol{P^{\rm DA}})}} \big[ \rho_s \eta \pi^{\rm gray} \big] -\sum_{\mathclap{s \in \mathcal{S}^{\rm RES}(\boldsymbol{P^{\rm DA}})}} \big[ \rho_s  \eta \pi^{\rm gray}  + \mu^{\rm green}_s \eta \pi^{\rm green} \big]  = 0 \\
    &\lambda^{\rm DA} + \mu^{\rm ub} - \eta \pi^{\rm gray}  \sum_{\mathclap{s \in \mathcal{S}}}  \rho_s  - \eta  \pi^{\rm green} \sum_{\mathclap{s \in \mathcal{S}^{\rm green}(\boldsymbol{P^{\rm DA}})}} \rho_s   - \eta \pi^{\rm green} \sum_{\mathclap{s \in \mathcal{S}^{\rm RES}(\boldsymbol{P^{\rm DA}})}}   \mu^{\rm green}_s   = 0 \\
    &\lambda^{\rm DA}  - - \eta \pi^{\rm gray}  - \eta  \pi^{\rm green} \sum_{\mathclap{s \in \mathcal{S}^{\rm green}(\boldsymbol{P^{\rm DA}})}} \rho_s   - \eta \pi^{\rm green} \sum_{\mathclap{s \in \mathcal{S}^{\rm RES}(\boldsymbol{P^{\rm DA}})}}   \mu^{\rm green}_s  = -\mu^{\rm ub} \leq 0 \\
    &\lambda^{\rm DA}  - \eta \pi^{\rm gray}  -\sum_{\mathclap{s \in \mathcal{S}^{\rm green}(\boldsymbol{P^{\rm DA}})}} \rho_s \eta  \pi^{\rm green}   \leq \sum_{\mathclap{s \in \mathcal{S}^{\rm RES}(\boldsymbol{P^{\rm DA}})}} \mu^{\rm green}_s \eta \pi^{\rm green} 
\end{align}
We derive the following bounds on the right-hand-side of this equation from the stationarity condition \eqref{stat_cond_2}:
\begin{align}
        0 \leq \eta \pi^{\rm green} \sum_{\mathclap{s \in \mathcal{S}^{\rm RES}(\boldsymbol{P^{\rm DA}})}}   \mu^{\rm green}_s  \leq \eta \pi^{\rm green}  \sum_{\mathclap{s \in \mathcal{S}^{\rm RES}(\boldsymbol{P^{\rm DA}})}}    \rho_s.
\end{align}
Given that $\boldsymbol{P^{\rm DA}} = P^{\rm h} \geq P^{\rm RES}_s, \ \forall s \in \mathcal{S}$, $\mathcal{S}^{\rm green}(\boldsymbol{P^{\rm DA}}) = \emptyset$, leading to the following upper-bound on the DA prices for which the day-ahead purchase quantity $\boldsymbol{P^{\rm DA}} = P^{\rm h}$ is an optimal solution:
\begin{align}
\lambda^{\rm DA}  \leq \eta \pi^{\rm gray} + \eta  \pi^{\rm green} \sum_{\mathclap{s \in \mathcal{S}^{\rm RES}(\boldsymbol{P^{\rm DA}})}}\rho_s,
\end{align}
that is,
\begin{align}
    \mathcal{P}(\boldsymbol{P^{\rm DA}}=P^{\rm h}) = \big] - \infty , \eta \pi^{\rm gray} + \eta  \pi^{\rm green} \sum_{\mathclap{s \in \mathcal{S}^{\rm RES}(\boldsymbol{P^{\rm DA}}=P^{\rm h})}}\rho_s \quad \ \ \big].
\end{align}
In addition, due to the definition of the unique and increasing scenarios, a scenario $s$ can only belong to $\mathcal{S}^{\rm RES}(\boldsymbol{P^{\rm DA}} = P^{\rm h})$ if $s=S$ and $P_S^{\rm RES} = P^{\rm h}$. In other words, these ranges can be defined explicitly as:
\begin{align}
   &  \mathcal{P}(\boldsymbol{P^{\rm DA}}=P^{\rm h}) = \big] - \infty , \eta \pi^{\rm gray} \big] , \text{ if } P_S^{\rm RES} < P^{\rm h} \\
   & \mathcal{P}(\boldsymbol{P^{\rm DA}}=P^{\rm h}) = \big] - \infty , \eta \pi^{\rm gray} + \eta  \pi^{\rm green} \rho_S\big] , \text{ if } P_S^{\rm RES} = P^{\rm h}
\end{align}

\end{subequations}

\newpage
\subsection{Derivation of the optimal bid curve} \label{appendix_bid_curve}

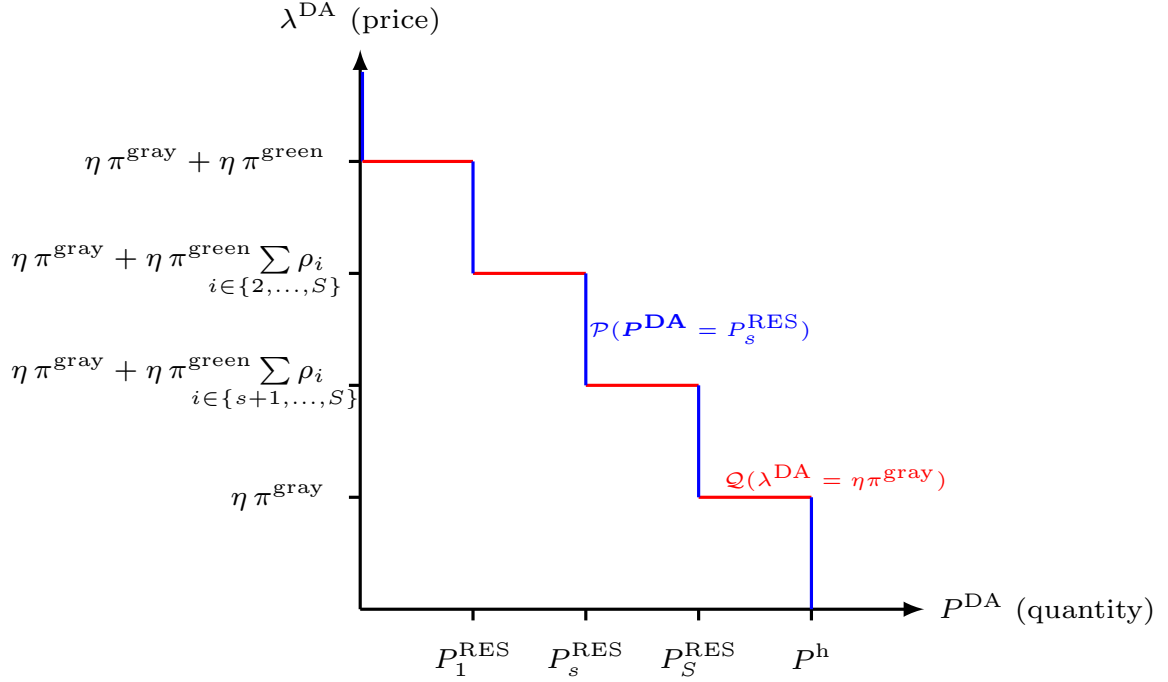
\begin{figure}[t]
  \centering
  \resizebox{\columnwidth}{!}{%
    \begin{tikzpicture}[
      thick,
      bidline/.style={blue},
      every node/.style={font=\scriptsize},
    ]

      \draw[->] (0,0) -- (5,0) node[anchor=west] {$P^{\rm DA}$ (quantity)};
      \draw[->] (0,0) -- (0,5) node[anchor=south] {$\lambda^{\rm DA}$ (price)};

      \foreach \x/\lbl in {
        1/P_{1}^{\rm RES},
        2/P_{s}^{\rm RES},
        3/P_{S}^{\rm RES},
        4/P^{\rm h}
      }
      {
        \draw (\x,0) -- (\x,-0.1) node[below=2pt] {$\lbl$};
      }

      \foreach \y/\lbl in {
        4/\eta\,\pi^{\rm gray} + \eta\,\pi^{\rm green},
        3/\eta\,\pi^{\rm gray} + \eta\,\pi^{\rm green}\sum\limits_{\mathclap{i\in\{2,\dots,S\}}}\rho_i,
        2/\eta\,\pi^{\rm gray} + \eta\,\pi^{\rm green}\sum\limits_{\mathclap{i\in\{s+1,\dots,S\}}}\rho_i,
        1/\eta\,\pi^{\rm gray}
      }
      {
        \draw (0,\y) -- (-0.1,\y) node[left=2pt] {$\lbl$};
      }

      \coordinate (P0) at (0.02,4);
      \coordinate (P1) at (1,4);
      \coordinate (P2) at (1,3);
      \coordinate (P3) at (2,3);
      \coordinate (P4) at (2,2);
      \coordinate (P5) at (3,2);
      \coordinate (P6) at (3,1);
      \coordinate (P7) at (4,1);
                \node[above=0pt of P7, yshift = -2pt, xshift = 6pt] {\color{red} \tiny $\mathcal{Q} ( \lambda^{\rm DA} = \eta \pi^{\rm gray} )$};
      \coordinate (P8) at (4,0);

      \draw[bidline] (0.02,4.8) -- (P0);

      \draw[red, thick]
         (P0) -- (P1) node[midway,below] {}
         (P2) -- (P3) node[midway,below] {}
         (P4) -- (P5) node[midway,below] {}
         (P6) -- (P7) node[midway,below] {};

      \draw[bidline]
         (P1) -- (P2) node[midway,right] {}
         (P3) -- (P4) node[midway,right,xshift=-3pt] {\tiny $\mathcal{P}(\boldsymbol{P^{\rm DA}}=P_s^{\rm RES})$}
         (P5) -- (P6) node[midway,right] {}
         (P7) -- (P8) node[midway,right] {};

    \end{tikzpicture}%
  }
  \caption{Step‐wise price–quantity bid curve.}
  \label{fig:step-curve-price-quantity-segments}
\end{figure}

To summarize, the conditions on the price-quantity pairs $(\lambda^{\rm DA},\boldsymbol{p^{\rm DA}})$ derived above describe the segments of the electrolyzer's optimal bid curve. In particular, the following ranges of day-ahead prices $\lambda^{\rm DA} \in \mathcal{P}(\boldsymbol{P^{\rm DA}})$ for which a given day-ahead purchase quantity $\boldsymbol{P^{\rm DA}}$ is optimal:
\begin{align}
    & \mathcal{P}(\boldsymbol{P^{\rm DA}}=0) = \begin{cases} 
    & \big[ \eta \pi^{\rm gray}   + \eta  \pi^{\rm green} (1 - \rho_1) , \ +\infty \big[ , \text{ if }P^{\rm RES}_1  = 0 \\
    & \big[ \eta \pi^{\rm gray}   + \eta  \pi^{\rm green}  , \ +\infty [ , \text{ if }P^{\rm RES}_1 >0.
    \end{cases} \\
    & \mathcal{P} (\boldsymbol{P^{\rm DA}}  = P_1^{\rm RES}) = \big[ \eta \pi^{\rm gray}  + \eta  \pi^{\rm green} \sum_{\mathclap{i \in \{2,...,S\}}} \rho_i \ , \eta \pi^{\rm gray} + \eta  \pi^{\rm green} \sum_{\mathclap{i \in \{1,...,S\}}} \rho_i \big] , \text{ if } P_1^{\rm RES} > 0 \\
    & \mathcal{P} (\boldsymbol{P^{\rm DA}}  = P_s^{\rm RES}) = \big[ \eta \pi^{\rm gray}  + \eta  \pi^{\rm green} \sum_{\mathclap{i \in \{s+1,...,S\}}} \rho_i \ , \eta \pi^{\rm gray} + \eta  \pi^{\rm green} \sum_{\mathclap{i \in \{s,...,S\}}} \rho_i \big] , \forall s \in \{2,...,S-1\} \\
    & \mathcal{P} (\boldsymbol{P^{\rm DA}}  = P_{S}^{\rm RES}) = \big[ \eta \pi^{\rm gray}  \ , \eta \pi^{\rm gray} + \eta \pi^{\rm green} \rho_S \big] , \text{ if } P_S^{\rm RES} < P^{\rm h} \\
    &  \mathcal{P}(\boldsymbol{P^{\rm DA}}=P^{\rm h}) = \begin{cases} 
    & \big] - \infty , \eta \pi^{\rm gray} \big] , \text{ if } P_S^{\rm RES} < P^{\rm h} \\
   & \big] - \infty , \eta \pi^{\rm gray} + \eta  \pi^{\rm green} \rho_S\big] , \text{ if } P_S^{\rm RES} = P^{\rm h}. 
   \end{cases}
\end{align}
each describe a \textit{vertical} segment of the stepwise price-quantity bid curve, for a fixed quantity $\boldsymbol{P^{\rm DA}}$. Similarly, the following ranges of optimal day-ahead purchase quantities $\boldsymbol{P^{\rm DA}}$ for a given day-ahead price $\lambda^{\rm DA}$:
\begin{align}
    & \mathcal{Q} ( \lambda^{\rm DA} = \eta \pi^{\rm gray}  + \eta  \pi^{\rm green} \sum_{\mathclap{i \in \mathcal{S}}} \rho_i ) = \big]0,P^{\rm RES}_1\big[,  \\
    & \mathcal{Q} ( \lambda^{\rm DA} = \eta \pi^{\rm gray}  + \eta  \pi^{\rm green} \sum_{\mathclap{i \in \{s,...S\} }} \rho_i ) = \big]P^{\rm RES}_{s-1},P^{\rm RES}_{s}\big[, \quad \forall s \in \{2,...,S\} \\
    & \mathcal{Q} ( \lambda^{\rm DA} = \eta \pi^{\rm gray} ) = \big]P^{\rm RES}_{S},P^{\rm h}\big[,
\end{align}
each describe a \textit{horizontal} segment of the stepwise price-quantity bid curve, for a fixed price $\lambda^{\rm DA}$. The resulting step-wise curve, from the price $\mathcal{P}$ and quantity $\mathcal{Q}$ segments, is illustrated in Figure \ref{fig:step-curve-price-quantity-segments}.

\end{document}